\documentclass[oneside,reqno,referee]{amsart}
\usepackage{amssymb}
\usepackage{verbatim}

\newtheorem{Th}{Theorem}[section]
\newtheorem{Cor}[Th]{Corollary}
\newtheorem{Lem}[Th]{Lemma}
\newtheorem{Prop}[Th]{Proposition}
\newtheorem{Claim}[Th]{Claim}

\newtheorem{Rem}[Th]{Remark}

\newtheorem{claim-num}{Claim}

\def\aut#1{\mbox{\rm Aut}(#1)}
\def\av#1{\overline{#1}}
\def\a{\alpha}

\def\cB{\mathcal B}
\def\cC{\mathcal C}
\def\cD{\mathcal D}
\def\cK{\mathcal K}
\def\cL{\mathcal L}
\def\cM{\mathcal M}
\def\cO{\mathcal O}
\def\cT{\mathcal T}
\def\cU{\mathcal U}
\def\cX{\mathcal X}
\def\cY{\mathcal Y}
\def\cZ{\mathcal Z}

\def\eps{\varepsilon}
\def\f{\varphi}
\def\N{\mathbf N}
\def\nc{\mbox{\rm NC}}
\def\Q{\mathbf Q}

\def\s{\sigma}

\def\str#1{\langle#1\rangle}
\def\stuple#1{#1_1,\ldots,#1_s}
\def\sym#1{\mbox{\rm Sym}(#1)}
\def\inn#1{\mbox{\rm Inn}(#1)}
\def\iat{\mbox{\rm IA}}
\def\ia#1{\mbox{\rm IA}(#1)}

\def\inv{{}^{-1}}
\def\sle{\subseteq}
\def\Mod#1{\ (\mbox{\rm mod}\ #1)}
\def\Z{\mathbf Z}
\def\id{\mbox{\rm id}}
\def\b{\beta}
\def\sge{\supseteq}

\def\text#1{\mbox{\rm #1}}

\renewcommand{\le}{\leqslant}
\renewcommand{\ge}{\geqslant}

\begin{document}

\title[Completeness of the automorphism groups]
{Infinitely generated free nilpotent groups: completeness of the automorphism groups}
\author{Vladimir Tolstykh}
\address{Vladimir Tolstykh\\ Department of Mathematics\\ Yeditepe University\\
34755 Kay\i\c sda\u g\i \\
Istanbul\\
Turkey}
\email{vtolstykh@yeditepe.edu.tr, tvlaa@rambler.ru}
\subjclass[2000]{20F28 (20F18, 03C60)}

\maketitle

\begin{abstract}
We transfer the results of Dyer, Formanek and Kassabov
on the automorphism towers of finitely generated free nilpotent
groups to infinitely generated free nilpotent groups.
We prove that the automorphism groups of infinitely
generated free nilpotent groups are complete.
By combining the results of Dyer, Formanek, Kassabov
with the results in the present paper, one gets that
the automorphism tower of any free nilpotent
group terminates after finitely many steps.
\end{abstract}

\maketitle


\section*{Introduction}

Baumslag conjectured in the 1970s that the automorphism
tower of a finitely generated free group (free nilpotent group)
must be very short. Dyer and Formanek \cite{DFo} justified
the conjecture concerning finitely generated free groups in the ``sharpest sense''
by proving that the automorphism group $\aut{F_n}$ of a
non-abelian free group $F_n$ of finite rank $n$ is complete. Recall that
a group $G$ is said to be {\it complete} if $G$ is centreless
and all automorphisms of $G$ are inner; it then follows
that $\aut G \cong G.$ Thus $\aut{\aut{F_n}} \cong \aut{F_n},$
or, in other words, the height of the automorphism tower
over $F_n$ is two. The proof of completeness of $\aut{F_n}$ given by Dyer and Formanek in \cite{DFo}
has later been followed by the proofs given by Formanek \cite{Fo}, by Khramtsov
\cite{Khr}, by Bridson and Vogtmann \cite{BrVo}, and by the author
\cite{To:Towers}. The proof given in \cite{To:Towers} works for
arbitrary non-abelian free groups; thus the automorphism groups of infinitely generated
free groups are also complete.

Let $F_{n,c}$ denote a free nilpotent group of finite
rank $n \ge 2$ and of nilpotency class $c \ge 2.$
In \cite{DFo2} Dyer and Formanek studied the automorphism
towers of free nilpotent groups $F_{n,2}$
of class two. They showed that the group $\aut{F_{n,2}}$
is complete provided that $n \ne 3.$ In the case
when $n=3$ the height of the automorphism tower
of $F_{n,2}$ is three. The main result of \cite{To:Camb}
states that the automorphism group of any infinitely
generated free nilpotent group of class two is complete.

In \cite{DFo_solv_like} Dyer and Formanek proved
completeness of the automorphism groups of groups
of the form $F_n/R'$ where $R$ is a characteristic
subgroup of $F_n$ which is contained in the commutator
subgroup $F_n'$ of $F_n$ and $F_n/R$ is residually torsion-free nilpotent.

In his Ph. D. thesis \cite{Kass} Kassabov
found an upper bound $u(n,c) \in \N$ for the height of the automorphism
tower of $F_{n,c}$ in terms of $n$ and $c,$ thereby
finally proving Baumslag's conjecture on
finitely generated free nilpotent groups. By analyzing
the function $u(n,c)$ one can conclude that if $c$ is
small compared to $n,$ then the height of
the automorphism tower of $F_{n,c}$ is at most three.

The main result of the present paper generalizing the main result
of \cite{To:Camb} states that the automorphism group
of any infinitely generated free nilpotent group
of nilpotency class $\ge 2$ is complete.

We would like, before discussing the structure of the
paper, to discuss model-theoretic terminology that will be used throughout
the paper. Generally
speaking, it is often required, when studying
the automorphisms of a given structure $\cM,$ to show
that a certain subset (relation) $S$ of $\cM$ is invariant under
all automorphisms of $\cM.$ Similarly,
it might be required to show that $S$ is fixed setwise
by all automorphisms of $\cM$ that fix a given subset
$R$ of $\cM$ pointwise. {\it Natural} examples
of subsets of $\cM$ that are invariant
under all automorphisms are given by subsets
that are definable in $\cM$ by formulae of some
logic $\cL.$ Let $\cL$ be a logic and $\cM$ a
structure in the language of $\cL.$ A subset $S$ of
$\cM$ is said to be {\it definable} in $\cM$ by means
of $\cL,$ if there is a formula $\chi(v)$ of $\cL$
such that $S=\chi(\cM),$ that is, $S$ is the set of
all realizations (``solutions'') of $\chi$ in $\cM;$
$S$ is said to be {\it definable with parameters} $\av
a$ where $\av a$ is a tuple of elements of $\cM$ if $S=\vartheta(\cM,\av a)$ for some formula $\vartheta(v,\av u)$
of $\cL.$

\begin{Prop} {\em \cite[Lemma 2.1.1]{Hodges}}
Suppose that a subset $S$ is definable
in $\cM$ {\em (}resp. definable in $\cM$ with
parameters $\av a${\em)} by means of some
logic $\cL.$ Then $S$ is fixed setwise
by all automorphisms of $\cM$ {\em(}resp. by all
automorphisms of $\cM$ that fix parameters
$\av a$ elementwise{\em).}
\end{Prop}

For instance, the family $i(G)$ of all involutions
is definable in any group $G$ by means of both first-order
and second-order logics (for the definition of the full second-order
logic see Section 2.8 of \cite{Hodges}); on the
other hand, the subgroup generated by $i(G)$
may not be first-order definable in $G,$
but it is second-order definable in $G$ etc.
In most of cases below, we prefer to
speak simply of definable subsets,
skipping references to logics whose formulae could be used to define
them, although it is safe to think that everywhere below ``definable'' means
``second-order definable'' (we must warn the reader that in model-theoretic texts
``definable'' often means ``first-order definable'').

We would like to stress that despite the presence of model-theoretic
terminology in the paper, {\it no} prior knowledge of model theory is actually assumed,
since the reader who is not familiar with
model theory can always substitute his or her own arguments to see that
a particular definable set behaves as in
the conclusion of the Proposition above.

The paper is organized as follows. Let $N$ be an infinitely generated free nilpotent
group of nilpotency class at least two and let
$\Delta \in \aut{\aut N}.$ In the first section
we show that the subgroup $\inn N$ of all
inner automorphisms of $N$ is definable
in the group $\aut N.$ This implies that
$\Delta$ can be followed by an inner automorphism $T_\rho$
of the group $\aut N$ so that the resulting
automorphism $\Delta_1=T_\rho \circ \Delta$
fixes the subgroup $\inn N$ pointwise.
The main result of
the second section states that the subgroup $\iat_2(N),$ the kernel
of the homomorphism $\aut N \to \aut{N/[N,N,N]},$
determined by the natural homomorphism $N \to N/[N,N,N],$
is also fixed pointwise by $\Delta_1;$ this constitutes
the major step towards the proof that $\Delta_1$
fixes all elements of the subgroup $\ia N$
of $\aut N.$ Sections 3 and 4 are devoted to the reconstruction of primitive
elements of $N$ in the group $\Gamma=\aut N,$ the key
result here says that the stabilizer
$$
\Gamma_{(x)}=\{\s\in \aut N : \s(x)=x\}
$$
where $x$ is a primitive element of $N$
is definable in $\aut N$ with a certain
parameter $\pi \in \aut N.$
In the final Section 5 we show that it is possible to reconstruct
in $\aut N$ the multiplication of ``independent'' primitive
elements. This enables us to show that the
automorphism $\Delta_1$ above
can be followed by an inner automorphism
of $\aut N$ so that the resulting
automorphism $\Delta_2$ fixes pointwise
some generating set of the group $\aut N.$

\section{Stabilizing conjugations}

Throughout the paper, $N$ will denote an infinitely
generated free nilpotent group and $c$ will always denote
nilpotency class of $N.$ We shall assume, if not
otherwise stated, that $c \ge 2.$

The aim of this section is to establish definability
of the subgroup $\inn N$ of all inner automorphisms of $N$ in
the group $\aut N.$ For convenience's sake, we shall simply call inner
automorphisms of $N$ conjugations. By the agreement,
if $a \in N$ then conjugation $\tau_a$ determined by $a$
is the automorphism
$$
\tau_a (z) =aza\inv \quad (z \in N).
$$

We use standard commutator notation: if
$\stuple a$ are elements of any group,
then $[a_1,a_2]=a_1a_2a_1^{-1} a_2^{-1}$
and $[a_1,\ldots,a_s] = [[a_1,\ldots,a_{s-1}],a_s]$
where $s \ge 3.$ Given a group $G,$ we denote
by $\gamma_k(G)$ the $k$-th term of the lower
central series of $G$: $\gamma_1(G)=G$ and
$\gamma_{k+1}(G)=[\gamma_k(G),G]$ $(k \in \N).$ We denote
the group $\gamma_k(N)$ by $N_k.$
Since $N \cong F/\gamma_{c+1}(F)$ where
$F$ is a free group having the same
rank as $N,$ the quotient
group $N/N_{k+1}$ where $1 \le k \le c-1$ is a free nilpotent
group of nilpotency class $k$; if $\varphi_k$ is the homomorphism
$\aut N \to \aut{N/N_{k+1}},$ determined by the natural
homomorphism $N \to N/N_{k+1},$ then the kernel
of $\varphi_k$ will be denoted by $\iat_k(N).$
Recall that the elements of the group $\iat_1(N)$
are called {\it IA-automorphisms} of $N$ and the standard
notation for the group $\iat_1(N)$ is $\ia N.$

In the following proposition we collect some
well-known facts about $N$ (its automorphism
group $\aut N$) that we shall frequently use below.

\begin{Prop} \label{BasicProps_o_N}
{\em (i)} Given any basis $\cX$ of $N$ and any family
$\{t_x : x \in \cX\}$ of elements of $N_{k}$ where $k \ge 2,$
there is an automorphism $\a \in \iat_{k-1}(N)$
that takes $x$ to $x t_x$ for every $x \in \cX;$

{\em (ii)} the $c$-th term $N_c$ of the lower
central series of $N,$ a free abelian
group, is the centre
of the group $N;$

{\em (iii)} any IA-automorphism of $N$
fixes each element of $N_c$ and any
element of $\iat_{c-1}(N)$ fixes each
element of $N',$ the commutator
subgroup of $N;$

{\em (iv)} the group $\iat_{c-1}(N)$ is the
centre of the group $\ia N,$
and then $\iat_{c-1}(N)$ is a {\em(}torsion-free{\em)} abelian group;

{\em (v)} the homomorphism $\aut N \to \aut{N/N_{k+1}}$
where $1 \le k \le c-1,$ determined by the natural
homomorphism $N \to N/N_{k+1},$ is surjective.
\end{Prop}

\begin{proof}
(i) See Theorem 31.25 in \cite{HN}.

(ii) See Theorems 31.61, 31.63 of \cite{HN}.

(iii) Both statements immediately follow from Theorem 5.3
of \cite{MKS}.

(iv) Let $\a$ be an arbitrary IA-automorphism
and $\eta$ an arbitrary element of $\iat_{c-1}(N).$
Take a basis $\cX$ of $N.$ Then $\a x =x t_x$
where $t_x \in N'$ and $\eta x =x s_x$
where $s_x \in N_c$ $(x \in \cX).$ By (iii), for
all $x \in \cX$ we have that $\eta\inv x =x s_x^{-1}.$
Again by (iii),
$$
\eta \alpha \eta\inv x = \eta \alpha (xs_x^{-1})=\eta( x t_x s_x^{-1})=x t_x =\alpha x
$$
for every $x \in \cX.$

Conversely, suppose that $\s$ is an element of the
centre of $\ia N.$ Then for every conjugation
$\tau_a \in \inn N \sle \ia N$ we have that
$$
\tau_a=\s \tau_a \s\inv=\tau_{\s(a)}
$$
whence, by (ii), $\s a \equiv a \Mod{N_c},$
or, in other words, $\s$ is an element of $\iat_{c-1}(N).$

(v) By (i), any basis of the free nilpotent
group $N/N_{k+1}$ can be lifted to a basis
of $N,$ whence the result.
\end{proof}

\begin{Prop} \label{Centreless}
The group $\aut N$ is centreless.
\end{Prop}

\begin{proof} We shall use the following claim in due course below;
it can be applied in the proof of the Proposition as well.

\begin{Claim} \label{First_app_o_Pi_star}
Let $\cX$ be any basis of $N$ and $\pi$
be an automorphism of $N$ such that $\pi$
preserves $\cX,$ fixes $x$ and all $\pi$-orbits on $\cY=\cX \setminus \{x\}$
are infinite. Suppose $\a$ is an $\iat$-automorphism
that commutes with $\pi.$ Then
$\a x=x.$
\end{Claim}

\begin{proof}
Write $\alpha x$ as $x w(x,\av y)$
where $w$ is a reduced word over $\cX,$ $w(x,\av y) \in N',$
and $\av y$ a finite tuple of elements of $\cY.$

Then the condition $\pi \alpha \pi\inv x =\alpha x$ implies that
$$
w(x, \av y)=\pi^k w(x, \av y)=w(x, \pi^k \av y).
$$
for all $k \in \Z.$ Since evidently there is a power $\pi^k$
of $\pi$ with $\av y \cap \pi^k \av y=\varnothing,$
then the word $w$ must be trivial.
\end{proof}

Now take any basis $\cX$ of $N$ and choose for every
$x \in \cX$ an automorphism $\pi_x$ of $N$ such that
$\pi_x(\cX)=\cX$ and all $\pi_x$-orbits on
$\cX$ but the orbit of $x$ which is equal to $\{x\}$ are infinite.

Suppose $\s$ is a central element of $\aut N.$
By Proposition \ref{BasicProps_o_N} (iv), $\s \in \iat_{c-1}(N).$
As $\s$ commutes with every $\pi_x,$ we have by
Claim \theClaim\ that $\s x=x$ for all $x \in \cX.$
\end{proof}

\begin{Th} \label{Def_o_conjs}
The subgroup $\inn N$ of all inner
automorphisms of $N$ is a definable
subgroup of the group $\aut N.$
\end{Th}

\begin{proof} We use induction on nilpotency
class $c$ of $N.$  When $c=2,$ the subgroup
$\inn N$ is first-order definable in the group $\aut N$
by Corollary 3.2 of \cite{To:Camb}.

The key tool for the induction step is provided by the following
result from \cite{To:ContMath}. Let $\widehat N$ denote
the free nilpotent group $N/N_c$ of nilpotency
class $c-1.$

\begin{Th} \label{Def-of-K} {\em \cite[Prop. 2.6]{To:ContMath}}
The group $\iat_{c-1}(N),$ the kernel
of the homomorphism $\widehat{\phantom m} : \aut N \to \aut{\widehat N},$
determined by the natural homomorphism $N \to \widehat N,$
is first-order definable in the group $\aut N.$
\end{Th}

\begin{Cor} \label{I-Is-Def}
Let $c \ge 3.$ Suppose that $\inn{\widehat N}$ is a definable
subgroup of $\aut{\widehat N}.$ Then any term $\gamma_k(\inn N)$ where $k \ge 2$ of the lower
central series of $\inn N$ is definable in the
group $\aut N.$
\end{Cor}

\begin{proof}
By Theorem \ref{Def-of-K} the preimage, say $L$ of
the group $\inn{\widehat N}$ under
$\widehat{\phantom a}$ is definable in the group $\aut N.$
Now $L=\inn N \cdot \iat_{c-1}(N).$ Since the elements of
$\iat_{c-1}(N)$ commute with all IA-automorphisms (Proposition \ref{BasicProps_o_N}),
the commutator subgroup $[L,L]$ of $L$
coincides with that one of $\inn N$:
$$
[L,L]=[\inn N,\inn N].
$$
We then apply induction on $k.$
\end{proof}

Recall that a {\it primitive} element
of a relatively free group $G$ is
one that can be included in some basis
of $G.$ We write $\nc(\s)$ for the normal closure
of a $\s \in \aut N$ in the group
$\aut N.$

\begin{Lem} \label{PrimConjs:Lifting}
Let $c \ge 3.$ An automorphism $\sigma$ of $N$ is conjugation by a
primitive element of $N$ if and only if

{\em (a)} $\widehat{\sigma}$ is conjugation
by a primitive element of $\widehat N;$

{\em (b)} the subgroup $\nc(\s)$ contains
no elements of the set $\iat_{c-1}(N) \setminus \gamma_{c-1}(\inn N).$
\end{Lem}

\begin{proof}
The necessity part is trivial. Let us prove the converse.
An argument similar to one used
in the proof of Corollary \ref{I-Is-Def} shows
that $\gamma_{c-1}(\inn N)$ is contained in $\nc(\sigma).$
Take an automorphism $\sigma$
that satisfies the conditions (a) and (b). Then
$\sigma$ has the form $\sigma =\tau_x \gamma$ for
some primitive $x \in N$ and some automorphism
$\gamma$ in $\iat_{c-1}(N).$
Suppose, towards a contradiction,
that $\sigma$ is not conjugation. Hence $\gamma \notin \gamma_{c-1}(\inn N).$

Observe that if $a$ is a primitive element of $N$ and $\tau_a \gamma_1$ and $\tau_a \gamma_2$
where $\gamma_k \in \iat_{c-1}(N)$ ($k=1,2$) are both in $\nc(\sigma),$ then
$\gamma_1 \equiv \gamma_2 \Mod{\gamma_{c-1}(\inn N)};$ indeed, otherwise the element
$
(\tau_a \gamma_1) (\tau_a \gamma_2)\inv
$
is a member of $\iat_{c-1}(N) \setminus \gamma_{c-1}(\inn N),$ contradicting (b).

Let $\cX$ be a basis of $N$ which contains $x.$
Write $\cY$ for $\cX \setminus \{x\}.$
Suppose that $\gamma x =x w(x,\av y)$ where $w$ is a
reduced word over $\cX$ and $\av y$ a finite tuple
of elements of $\cY.$ Consider
an automorphism $\pi$ of $N$ which preserves
$\cX,$ fixes $x$ and such that tuples $\pi \av y$
and  $\av y$ have no common elements. As $\pi x =x$
we have by the observation above that
$
\gamma^\pi =\gamma \tau_s
$
where $s \in N_{c-1}.$ It follows that
\begin{equation}
x w(x, \pi \av y) =sxs\inv w(x,\av y).
\end{equation}
Consider an endomorphism
$\eps$ of $N$ which takes all elements of $\pi \av y$
to $1$ and fixes all other elements of $\cX.$ By applying
$\eps$ to the both parts of (\theequation), we get
that $w(x,\av y) = tx\inv t\inv x$ where $t =\eps(s).$
As $t x\inv t\inv x =[t, x\inv]=[t\inv,x]$ is
in the centre of $N,$ we have
that
$$
\gamma x =x w(x,\av y)=x (t\inv x t x\inv) = (t\inv x t x\inv) x =
t\inv x t.
$$
Hence the automorphism $\tau_t \gamma$ fixes $x.$
Now $\s = \tau_x \gamma = \tau_{xt\inv } \tau_t \gamma,$
and conjugating $\s$ by any IA-automorphism
that sends $xt\inv$ to $x,$ we obtain that
$\tau_x (\tau_t \gamma)$ is in $\nc(\sigma).$
Thus, without loss of generality, we may assume
that $\gamma x=x.$

It is easy to see that $x$ (any primitive element
of $N$) commutes with an element $s \in N_{c-1}$
if and only if $s \in N_c,$ the centre
of $N.$ Let $\rho$ be any automorphism of $N$ which
takes $x$ to itself. Again, $\gamma^\rho = \gamma \tau_s$
for a suitable $s \in N_{c-1},$ whence $x =sxs\inv,$ and
then $s \in N_c.$ Therefore
$\tau_s=\id,$ or, in other words, $\gamma$ commutes with $\rho.$

Consider an arbitrary element $y \in \cY.$
Assume that $\gamma y = y w_y(y,x,\av z)$ where $w_y$
is a reduced word over $\cX$ and $\av z$
is a finite tuple of elements of $\cX \setminus \{x,y\}.$
As it is easy to find an automorphism $\pi$ of $N,$
preserving $\cX,$ fixing both $x$ and $y$ (and
hence commuting with $\gamma$), and such that
$\pi \av z \cap \av z=\varnothing,$ the tuple $\av z$
must be empty. Further, by considering an automorphism
interchanging distinct elements $y_1,y_2 \in \cY$ and
fixing all elements in $\cX \setminus \{y_1,y_2\},$
we see that there is a group word $w(*_1,,*_2)$
over an alphabet having empty intersection
with $N$ such that $w(y,x)=[w(*_1,*_2)]^{*_1,}_{y,}{}^{*_2}_x=w_y(y,x)$ for every $y \in \cY$
(or, somewhat more formally, there is a term $w(*_1,*_2)$
in two variables of the language of group theory
such that its value $w(y,x)$ at $(y,x)$ is $w_y$ for every $y \in \cY$).

Also, taking distinct $y,z \in \cY$ and considering
the automorphism
$\rho$ of $N$ with $\rho y=yz,$
fixing $\cX \setminus \{y\}$ pointwise,
 we obtain that
\begin{equation} \label{firstcoordlin}
w(yz,x)=w(y,x) w(z,x).
\end{equation}

Fix a $y \in \cY.$ Consider the automorphism $\pi$
of $N$ interchanging $x$ and $y$ and fixing pointwise
$\cX \setminus \{x,y\},$ the automorphism $\rho$
which takes $x$ to $xy$ and fixes $\cY$ pointwise
and any IA-automorphism $\a$ taking $xyx\inv$ to $y.$
We then have that the automorphism
$$
[(\tau_x \gamma)^\rho (\tau_x\gamma)\inv]^\alpha =
[\tau_{xyx\inv} \gamma^\rho \gamma\inv]^\alpha = \tau_y \gamma^\rho \gamma\inv
$$
belongs to $\nc(\s),$ and so does the automorphism $\tau_ y \gamma^\rho =(\tau_x \gamma)^\pi.$
Therefore for some $s \in N_{c-1}$
\begin{equation}
\gamma^\rho \gamma\inv = \tau_s \gamma^\pi.
\end{equation}
One readily verifies that (\theequation) implies
that
\begin{equation}
z w(z,xy) = s z s\inv w(z,x) w(z,y)
\end{equation}
for every $z \in \cY.$ In particular,
$$
y w(y,xy) = s y s\inv w(y,x).
$$
Let $\eps$ be the endomorphism of $N$ which
sends $x$ to $y$ and fixes $\cY$ pointwise. Then
$y = \eps(s) y \eps(s)\inv,$ and hence $\eps(s) \in N_c.$
Fix a $z \in \cY \setminus\{y\}$ and apply
$\eps$ to the both parts of (\theequation), thereby
getting that
$$
z w(z,y^2)=z w(z,y)^2, \text{ or } w(z,y^2) =w(z,y)^2.
$$
By (\theequation), we also have that $w(z^2,y)=w(z,y)^2.$
We then claim that $w(z,y)=1.$
We define a function $\nu_z$
on the set of all basic commutators over $\{z,y\}$
taking values in $\N$:
$\nu_z(z)=1, \nu_z(y)=0$ and $\nu_z([c_1,c_2])=
\nu_z(c_1) +\nu_z(c_2)$ where $c_1,c_2$
are basic commutators over $\{z,y\}.$
The function $\nu_y$ is defined similarly.
Let $b=b(z,y)$ be a basic commutator
of weight $c$ over $\{z,y\}.$ Observe that
$\nu_z(b)+\nu_y(b)=c$  and
$$
b(z^k,y)b(z,y^k)=b(z,y)^{k^{\nu_z(b)}}
b(z,y)^{k^{\nu_z(b)}}=b(z,y)^{k^c}
$$
where $k$ is a natural number.
Therefore
\begin{equation} \label{glue_em_to_kill_em}
v(z^k,y)v(z,y^k)=v(z,y)^{k^c} \qquad (k \in \N)
\end{equation}
for every element/word $v \in \str{z,y} \cap N_c.$
It then follows that
$$
w(z,y)^4=w(z,y)^2 w(z,y)^2=w(z^2,y)w(z,y^2)=w(z,y)^{2^c}
$$
As $4 < 2^c,$ $w(z,y)=1.$ Then $\gamma =\id,$
which completes the proof of the lemma.
\end{proof}

As it was said above, the group $\inn N$ is definable
in $\aut N$ when $c=2;$ then so is the family of all conjugations
by primitive elements, for $\inn N=\nc(\tau)$ where $\tau \in \inn N$
if and only if $\tau$ is conjugation by a primitive element.
Then by Lemma \ref{PrimConjs:Lifting} conjugations by primitive
elements are definable in the group $\aut N$
for any $c \ge 2,$ and hence the subgroup $\inn N$
is definable in the group $\aut N$ for any $c \ge 2.$
\end{proof} 

Now we obtain some corollaries of Theorem \ref{Def_o_conjs}.

\begin{Cor}
All subgroups $\iat_k(N)$ where $1 \le k \le c-1$
are definable in the group $\aut N.$
\end{Cor}

\begin{proof}
By Theorem \ref{Def-of-K} and by Theorem \ref{Def_o_conjs}.
\end{proof}

\begin{Prop} \label{StabilizingConjugations}
{\em (i)} Let $\Delta \in \aut{\aut N}.$ Then there is an inner automorphism $T_\pi
\in \aut{\aut N}$ such that the composition $T_\pi \circ \Delta$
fixes the subgroup $\inn N$ pointwise;

{\em (ii)} suppose that $\Delta \in \aut{\aut N}$
fixes every element of $\inn N.$ Then
$$
\Delta \sigma \equiv \sigma \Mod{\iat_{c-1}(N)}
$$
for all $\s \in \aut N.$
\end{Prop}

\begin{proof}
(i) We saw above that conjugations
by primitive elements of $N$ are definable;
a generating set $X$ of $\inn N$ consisting
of conjugations by primitive elements is a basis of $\inn N$ if and only
if
$$
\tau \not\in \str{X \setminus \{\tau\}}\quad \forall \tau \in X
$$
where $\str{X \setminus \{\tau\}}$ is the subgroup
of $\aut N$ generated by the set $X \setminus \{\tau\}.$
Thus every automorphism of the group $\aut N$ takes an
arbitrary basis of $\inn N$
to another basis of $\inn N.$

Fix a basis $X$ of $\inn N$ and a $\Delta \in \aut{\aut N}.$
As $\Delta(X)$ is a basis of $\inn N,$
there is a $\rho \in \aut N$ such that
$$
\Delta(\tau) =\rho \tau \rho\inv \quad (\tau \in X).
$$
Now if $T_{\rho\inv}$ is the inner automorphism of the group $\aut N$ determined
by $\rho\inv,$ that is, if
$$
T_{\rho\inv}(\s) = \rho\inv \s \rho \quad (\s \in \aut N),
$$
then $T_{\rho\inv} \circ \Delta$ fixes $X$
pointwise.

(ii). Let $z$ be any element of $N$ and $\s \in \aut N.$ We have that
$$
\tau_{\sigma z} =\Delta(\tau_{\sigma z})=\Delta(\sigma \tau_z \sigma\inv) =
\Delta(\sigma) \tau_z \Delta(\sigma)\inv=\tau_{\Delta(\s) z}.
$$
Then $\s z \equiv \Delta(\s) z \Mod{N_c}$ for every $z \in N,$
and hence $\s\inv \Delta(\s) \in \iat_{c-1}(N).$
\end{proof}

\section{Stabilizing IA-automorphisms}

Let $\cX$ be any basis of $N.$ Then by $\sym\cX$ (as it is,
for instance, done in \cite{BrRo}) we denote the group
of all automorphisms of $N$ that preserve $\cX$ setwise.
For simplicity's sake, we call elements of $\sym\cX$
{\it permutational} automorphisms relative to $\cX.$
Similarly, if $X$ is a basis of the group $\inn N$ (a basis
set of conjugations), the symbol $\sym X$ denotes the group of all automorphisms
of $N$ acting on $X$ as permutations.

\begin{Rem} \label{Chto_te_chto_eti}
\em Suppose that the elements of a basis $\cX$
of $N$ determine the elements of a given basis
set of conjugations $X$: $X=\{\tau_x : x \in \cX\}$
(in similar situations below we shall also say
that a basis of $N$ {\it determines} a given
basis of $\inn N$). Then for any $\rho \in \sym X$ there is a uniquely
determined $\sigma \in \sym \cX$ such that the actions of
$\rho$ and $\sigma$ on $\ia N$ are the same:
\begin{equation}
\rho \alpha \rho\inv = \sigma \alpha \sigma\inv \quad
\end{equation}
for every $\alpha \in \ia N.$ Indeed, it follows from the assumptions made that there
is an $\eta \in \iat_{c-1}(N)$ with $\rho =\s \eta.$ Then
(\theequation) is a corollary of the fact that elements
of $\iat_{c-1}(N)$ commute with all elements of $\ia N$ (Proposition
\ref{BasicProps_o_N}).

In particular, every time when we need an element of $\sym X$ that acts
on $\ia N,$ we can use a suitable element
of $\sym\cX$ and vice versa.
\end{Rem}

We use below the notation from the theory of permutation
groups. Recall that given a group $G$ acting on a set
$A$ we denote by $G_{(B)}$ and $G_{\{B\}}$ the
pointwise and the setwise stabilizers of a subset
$B$ of $A,$ respectively. Any symbol of the form
$G_{*_1,*_2}$ is the intersection of subgroups
$G_{*_1}$ and $G_{*_2}.$ It is convenient,
when working with subgroups $G$ of $\aut N$
and when using some basis $\cB$ of $N,$ to denote simply by $G_{\str{\cC}}$
where $\cC \sle \cB$ the subgroup $G_{(\cB \setminus \cC),\{\str{\cC}\}}$
(consisting of elements of $G$ that fix the set
$\cB \setminus \cC$ pointwise and preserve
the subgroup of $N$ generated by $\cC).$

\begin{Prop} \label{IA_k-stabs_Are_Def}
Let $x$ be a primitive element of $N$
and let $k$ with $1 \le k \le c-1$ be a natural number.
Then the group $\iat_k(N)_{(x)},$ consisting
of the elements of the group $\iat_k(N)$ stabilizing
$x,$ is definable with the parameter $\tau_x$ in the group $\aut N$
\end{Prop}

\begin{proof} Let $X$ be any basis set of conjugation
of which $\tau_x$ forms a part. Write $Y$ for
$X \setminus \{\tau_x\}.$ Choose also a basis
$\cX$ of $N$ such that the elements of $\cX$
determine the elements of $X$; we assume that $x \in \cX.$

Consider an element $\pi^*$ of $\sym X$ which, under
the conjugation action, fixes $\tau_x$ and
moves elements of $Y$ as a permutation having no fixed points
and having infinitely many infinite orbits. According to Remark \ref{Chto_te_chto_eti},
we may assume without loss of generality (for we are going
to consider conjugation actions of elements of $\aut N$
on $\iat(N)$) that $\pi^* x=x$ and that infinitely
many orbits of the action of $\pi^*$
on $\cY= \cX \setminus \{x\}$ are all infinite.

By Claim \ref{First_app_o_Pi_star}, if an $\eta \in \iat_{c-1}(N)$ commutes with $\pi^*$
(with any element of the coset $\pi^* \iat_{c-1}(N)$), then $\eta x=x.$

When using symbols like $G_{(*)}$ with $G=\sym\cX$
or $G=\sym X,$ we shall replace $\sym \cX$
and $\sym X$ by $\Pi$ and $P,$ respectively.
By $Z^+(\pi^*)$ we will denote the subgroup of the centralizer
$Z(\pi^*)$ of $\pi^*$ in the group $\aut N$ consisting of automorphisms
sending $x$ to itself.

Let $L$ be a free nilpotent
group, $\cB$ a basis of $L$ and
$\cB = \bigcup_{i \in I} \cB_i$ be a partition of
$\cB.$ If $\rho_i$ is an automorphism of the subgroup
$\str{\cB_i}$ ($i \in I$) we shall denote by
$\circledast_{i \in I} \rho_i$ the only automorphism
of $L$ which extends all automorphisms $\rho_i.$

The following result will be used in
Section \ref{Rec_o_PrimEls} as well. Recall that a {\it moiety} of a given infinite set $I$
is any subset $J$ of $I$ with $|J|=|I\setminus J|.$

\begin{Lem} \label{Populating_cZ_0}
Let $\cY_0 \sle \cY$ be a subset of $\cY$ maximal with
the property ``there are no elements lying in the same orbit
of $\pi^*$'' and let $\cZ_0$ be a moiety of $\cY_0.$
Take an arbitrary automorphism $\gamma$ of the
group $\str{x,\cZ_0}$ that fix $x.$ Then
there is a $\sigma \in Z^+(\pi^*)$ and
a $\rho \in \Pi_{(x)}=\sym \cX_{(x)}$ such that the
automorphism $\s \s^\rho$ fixes all elements of $\cX \setminus \cZ_0$
and coincides with $\gamma$ on $\cZ_0.$
In the case when $\gamma \in \iat_k(\str{x,\cZ_0}),$
$\s$ can be taken from $\iat_k(N).$
\end{Lem}

\begin{proof}
Write $\cY_k$ for $(\pi^*)^k \cY_0$ where $k \in \Z.$ Hence
$$
\cY =\bigcup_{k \in \Z} \cY_k
$$
is a partition of $\cY$ into moieties.

We then partition $\cY_0 \setminus \cZ_0$ into
countably many moieties indexed by elements
of $\N \setminus \{0\},$ getting that
$$
\cY_0 = \bigcup_{n \in \N} \cZ_n.
$$
Assume that $\cZ_k=\{z_{i,k} : i \in I\}$ where $k$ runs over
$\N.$

First, to illustrate the idea we consider the case when $\gamma$ preserves
the group $\str{\cZ_0};$ let $\nu$ denote
the restriction of $\gamma$ on $\str{\cZ_0}.$
Our goal is then to construct an automorphism $\sigma_0$
of $\str{\cY_0}$ that can be written as a $\circledast$-product of automorphisms
of the groups $\str{\cZ_k}$ of the form
$$
\underbrace{\nu \circledast \nu\inv \circledast \id} \circledast \underbrace{\nu \circledast \nu\inv \circledast \id} \circledast \ldots \circledast
\underbrace{\nu* \nu\inv \circledast \id} \circledast \ldots
$$
where the symbols $\nu,\nu\inv$ represent the action (up to an isomorphism
of actions) of $\sigma_0$ on subgroups $\str{\cZ_k}$ ($k \in \N$). Formally,
for every $k \in \Z$ we denote the bijection
$$
z_{i,0} \mapsto z_{i,k} \quad (i \in I)
$$
from $\cZ_0$ to $\cZ_k$ by $\lambda_k$ and then extend $\lambda_k$ to the isomorphism $\lambda^*_k$ of groups $\str{\cZ_0}$
and $\str{\cZ_k}.$ Then $\sigma_0$ is defined as follows:
\begin{equation} \label{def_o_eta0}
\begin{array}{ll}
\sigma_0(z_{i,k})= \lambda^*_k\gamma(z_{i,0}),         &\quad \text{ if } k \equiv 0 \Mod 3,\\
\sigma_0(z_{i,k})= \lambda^*_k\gamma\inv(z_{i,0}),     &\quad \text{ if } k \equiv 1 \Mod 3, \\
\sigma_0(z_{i,k})= z_{i,k}                             &\quad \text{ if } k \equiv 2 \Mod 3
\end{array}
\end{equation}
for every $i \in I.$

Using bijections $\lambda_k,$ one can construct permutational
automorphisms $\rho_0$ and $\rho_1$ of the group $\str{\cY_0}$
such that the conjugates $\sigma_0^{\rho_0}$ and $\sigma_0^{\rho_1}$ are
automorphisms of the form
$$
\underbrace{\nu\inv \circledast \nu \circledast \id} \circledast \underbrace{\nu\inv \circledast \nu \circledast \id} \circledast \ldots \circledast
\underbrace{\nu\inv \circledast \nu \circledast \id} \circledast \ldots
$$
and of the form
$$
{\id \circledast \nu \circledast \id} \circledast \underbrace{\nu\inv \circledast \nu \circledast \id} \circledast \ldots \circledast
\underbrace{\nu\inv \circledast \nu \circledast \id} \circledast \ldots
$$
respectively, and the products $\s_0^{\rho_0}\s_0$ and
$\s_0^{\rho_1}\s_0$ are the identity automorphism
and the automorphism
\begin{equation}
\nu \circledast \id \circledast \id \circledast \ldots \circledast \id \circledast \id \circledast \id \circledast \ldots
\end{equation}
of $\str{\cY_0},$ respectively.

It is then easy to extend $\sigma_0$ to the sets $\cY_k$ with $k \ne 0$
aiming at obtaining an element $\sigma$ of $Z^+(\pi^*)$ that
fixes $x$ and commutes with $\pi^*.$

There is a permutational automorphism $\rho \in \Pi_{(x)}$ such that
\begin{itemize}
\item $\rho$ fixes $x;$
\item the action of $\rho$ on any set $\cY_k$ where $k \ne 0$ is isomorphic
to that one of $\rho_0$ on $\cY_0;$
\item the restriction of $\rho$ on $\cY_0$ is $\rho_1;$
\item the product $\sigma \sigma^\rho$ acts trivially
on any $\cY_k$ where $k \ne 0$ and the restriction of $\rho$
on $\str{\cY_0}$ equals the automorphism (\theequation).
\end{itemize}
In particular, $\sigma \sigma^\rho \in \Gamma_{(\cX \setminus \cZ_0),\{\str{\cZ_0}\}}=
\Gamma_{\str{\cZ_0}}$ where $\Gamma=\aut N.$

Now it is easy to check that permutational
automorphisms $\rho_0,\rho_1$ relative
to $\cY_0$ and a permutational automorphism
$\rho \in \Pi_{(x)}$ can also be found in
the general case when the images $\gamma(z_{i,0}) \in \str{x,\cZ_0}$
of elements of $\cZ_0$ are not necessarily in $\str{\cZ_0}.$
\end{proof}

Assume some basis $\cU$ of $N_c$ consisting of basic commutators
of weight $c$ over $\cX$ is fixed.

Write $K$ for the group $\iat_{c-1}(N).$ Lemma \theLem\ implies that there is a moiety $\cZ_0$
of $\cY$ such that the subgroup $K_{\str{\cZ_0}}$ (the notation
here and up to the end of the proof is relative to the basis $\cX,$
if not otherwise stated) and the set, say $B_x$ of elements of $K_{(\cX \setminus \cZ_0)}$
which take $z_{i,0} \in \cZ_0$ either
to itself, or to an element $z_{i,0} t_i$
where all elements of $\cU$ participating
in the decomposition of $t_i$ have occurrences
of $x$ ($i\in I$), are both contained in the set
$$
\{\eta \eta^\rho : \eta \in \iat_{c-1}(N), \rho \in \Pi_{(x)}, [\eta,\pi^*]=\id\}
$$
which is equal to the set
$$
\{\eta \eta^\rho : \eta \in \iat_{c-1}(N), \rho \in P_{(\tau_x)}, [\eta,\pi^*]=\id\}
$$
and hence definable with the parameters $\tau_x,X$ and $\pi^*.$

The following statement is essentially
a part of the proof of Theorem 2.5 in \cite{To:JLM_Bergman}.

\begin{Claim}
Given a moiety $\cC$ of a basis $\cB$ of $N,$
every element of $K$ can be written
as a product of at most $2(c+1)$ conjugates
of elements of $K_{\str{\cC}}$
by permutational automorphisms relative
to $\cB.$
\end{Claim}

For the reader's
convenience, we reproduce the proof, for we shall need
below both the statement and the method of its proof to
be applied in somewhat different situations.

\begin{proof}
We can write $K$ as $K=K_{(\cB \setminus \cC)} \cdot K_{(\cC)}$ (Proposition
\ref{BasicProps_o_N} (i) and (iii)).
It suffices to prove that every element of $K_{(\cC)}$
is a product of at most $c+1$ conjugates of $K_{\str{\cC}}$
by permutational automorphisms, because the groups
$K_{(\cB \setminus \cC)}$ and $K_{(\cC)}$ are conjugate
by some permutational automorphism. Write $\cD$
for $\cB \setminus \cC$ and let $\cD=\{d_i : i \in I\}.$
Take an $\alpha \in K_{(\cC)}.$ Then $\a$ acts trivially
on $\cC$ and
$$
\a d_i =d_i s_i \qquad (i \in I)
$$
where $s_i \in N_c.$ We make use of
the fact that any basic commutator of weight $c$ over $\cB$
is formed from at most $c$ elements of $\cB.$ Thus
we write $\cC$ as a union of $c+1$ disjoint moieties:
$$
\cC = \cC_1 \cup \ldots \cup \cC_{c+1}.
$$
We write $s_i$ as a product of elements of some basis
of $N_c,$ basic commutators over $\cB.$ We then rewrite $s_i$ $(i \in I)$
in the form
$$
s_i =s_{i,1} \ldots s_{i,c+1}
$$
collecting into $s_{i,k}$ all those basis elements in
the decomposition of $s_i$ that do {\it not} contain,
as basic commutators, occurrences of letters from
$\cC_k$ ($k=1,\ldots,c+1).$ Then we define $c+1$
elements of $\iat_{c-1}(N)$ such that their product $\a_1 \ldots \a_{c+1}$ is equal
to $\a$:
$$
\begin{array}{lll}
\a_k d_i &=d_i s_{i,k}, &\quad (i \in I), \\
\a_k d   &=d,           &\quad (d \in \cC).
\end{array}
$$
As $\a_k \in K_{\str{\cD \cup (\cC \setminus \cC_k)}}=
K_{\str{\cB \setminus \cC_k}}$ $(k=1,\ldots,c+1),$
we have the desired, for the subgroups $K_{\str{\cB \setminus \cC_k}}$
and $K_{\str{\cC}}$ are conjugate by
a permutational automorphism.
\end{proof}

It follows from the Claim that the group $K_{\str{\cY}}$ is contained
in the group generated by conjugates of elements
of $K_{\str{\cZ_0}}$
by elements of $P_{(\tau_x)}.$ Similarly,
the proof of the Claim can be applied to show
that the subgroup generated by conjugates of
elements of $B_x \cup K_{\str{\cZ_0}}$ by elements of $P_{(\tau_x)}$
contains the family $C_x$ of all elements $\{\alpha\}$
of $\iat_{c-1}(N)$ such that $\alpha x=x$ and the image
any element $y \in \cY$ is of the form $yt_y$
where $t_y \in N_c$ is either $1,$ or again every
basis element of $\cU$ appearing in the decomposition
of $t_y$ has occurrences of $x.$ We complete the proof
for the group $K_{(x)}$ by observing that
$$
K_{(x)} =K_{\str{\cY}} \cdot C_x.
$$

Before going on, let us summarize the proof of definability
of $K_{(x)}=\iat_{c-1}(N)_{(x)}$ over $\tau_x$: we have proved
that there is a basis set $X$ of conjugations with $\tau_x \in X$, an element
$\pi^*$ of $P=\sym X$ acting on $X$ in the manner described
above such that the group $\iat_{c-1}(N)_{(x)}$ equals the subgroup
generated by conjugates of elements of $\iat_{c-1}(N)$ commuting
with $\pi^*$ by elements of $P_{(\tau_x)}$:
$$
\iat_{c-1}(N)_{(x)} =\str{\rho \eta \rho\inv : \rho \in P_{(\tau_x)}, \eta \in \iat_{c-1}(N), [\eta,\pi^*]=\id}.
$$

\begin{Lem} \label{I_k_in_the_LHS}
For every $k$ with $1 \le k \le c-1$
\begin{equation}
\iat_{k}(N)_{(x)} =\str{\rho \eta \rho\inv : \rho \in P_{(\tau_x)},
\eta \in \iat_{k}(N), [\eta,\pi^*]=\id}.
\end{equation}
\end{Lem}

\begin{proof}
We prove the statement by induction on $k.$
Suppose that (\theequation) is true for the group
$\iat_{k}(N)$ and $k \ge 2.$ Write $J$ for the
group
$$
\str{\rho \eta \rho\inv : \rho \in P_{(\tau_x)},
\eta \in \iat_{k-1}(N), [\eta,\pi^*]=\id}.
$$
Then $J \sle \iat_{k-1}(N)_{(x)}.$ It is also
clear that $\iat_{k}(N)_{(x)}$ is contained
in $J$ by the induction hypothesis.

We consider the natural homomorphism $\widetilde{\phantom m}$
from the group $N$ onto the free nilpotent group $\widetilde N =N/N_{k+1}$
of nilpotency class $k.$  The homomorphism
$\widetilde{\phantom m}$ determines the homomorphism
$\aut N \to \aut{\widetilde N}$ which we will denote
by the same symbol $\widetilde{\phantom m}.$
By the above considerations on the group $\iat_{c-1}(N)$
we have that
\begin{equation}
\iat_{k-1}(\widetilde N)_{(\widetilde x)} =
\str{\rho \eta \rho\inv : \rho \in P_{(\widetilde \tau_x)},
\eta \in \iat_{k-1}(\widetilde N), [\eta,\widetilde{\pi^*}]=\id}.
\end{equation}
Note that since the homomorphism $\aut N \to \aut{\widetilde N}$
is surjective (Proposition \ref{BasicProps_o_N} (iv)), then
$$
\widetilde J = \iat_{k-1}(\widetilde N)_{(\widetilde x)}.
$$

Let $\a \in \iat_{k-1}(N)_{(x)}.$ As $\a x=x,$ then $\widetilde \a (\widetilde x) =\widetilde x,$
and by (\theequation) $\widetilde \a \in \widetilde J.$ Hence
there exist an $\eta \in J$ and a $\beta \in \iat_{k}(N)$ such that
$\a =\beta \eta.$ Suppose that $\beta x =xt$ where $t \in N_{k+1}.$
We have that
$
\a x =\beta\eta(x)=\beta(x) =xt,
$
whence $t=1.$ Therefore $\beta \in \iat_{k}(N)_{(x)} \sle J,$
and we are done.
\end{proof}

The proof of Proposition \ref{IA_k-stabs_Are_Def} is now completed.
\end{proof}

We are going to make use of the following result from paper \cite{BrGu}
by Bryant and Gupta.

\begin{Th} \label{BrGu_T,deltas}
Let $F_{n,c}$ be a free nilpotent
group of rank $n \ge 3$ and of class $c \ge 3$ where $n \ge c-1,(c+1)/2.$
Let $x_1,\ldots,x_n$ be a basis of $F_{n,c}.$ Consider the automorphism
$\delta_c \in \aut{F_c}$ defined via
$$
\begin{array}{lll}
\delta_c: & x_1 &\mapsto x_1 [x_1,x_2,\ldots,x_{c-1},x_1],\\
		& x_i &\mapsto x_i, \quad i >1.
\end{array}
$$
Then the group $I_{c-1}=\iat_{c-1}(F_{n,c})$ is generated
by conjugates of elements of the set $(I'~\cap~I_{c-1})~\cup~\{\delta_c\},$ where $I=\iat(F_{n,c})$
and $I'$ is the commutator subgroup of $I,$ by tame automorphisms of $F_{n,c}.$
\end{Th}

\noindent (See \cite[p. 313, p. 317]{BrGu} for the details; for applications below, it will be enough to know that
$\iat_{c-1}(F_{n,c})$
is contained in the normal closure of the set $(I' \cap I_{c-1}) \cup \{\delta_c\}).$

In the first part of the proof of our next result,
we shall check that any automorphism $\Delta$ of the group $\aut N$ stabilizing
all conjugations, stabilizes an analogue of $\delta_c$
in the group $\aut N.$

\begin{Prop} \label{StabilizingIA3}
Let $\Delta$ be an automorphism of the group
$\aut N$ stabilizing all elements of the
group $\inn N.$ Then $\Delta$ fixes
the group $\iat_2(N)$ elementwise.
\end{Prop}

\begin{proof}
First, we observe that any automorphism $\Delta$ satisfying
the conditions of the Proposition takes any commutator
of elements of $\iat(N),$ and hence
any element of the commutator subgroup
of $\ia N,$ to itself. Indeed, let $\alpha_1,\alpha_2
\in \ia N.$ By Proposition \ref{StabilizingConjugations} (ii),
$\Delta(\alpha_k)=\alpha_k \eta_k$ where $\eta_k \in \iat_{c-1}(N),$
$k=1,2.$ Therefore
$$
\Delta(\alpha_1 \alpha_2 \alpha_1^{-1} \alpha_2^{-1})=
\alpha_1 \eta_1 \cdot \alpha_2 \eta_2 \cdot \eta_1^{-1} \alpha_1^{-1} \cdot \eta_2 ^{-1} \alpha_2^{-1}
=\alpha_1  \alpha_2 \alpha_1^{-1} \alpha_2^{-1},
$$
since $\eta_1,\eta_2$ are central elements of the group
$\ia N.$ In a similar vein, one proves that if
$\alpha \in \ia N$ is stabilized by $\Delta,$ then so is
any conjugate of $\alpha$ by an element of $\aut N.$

Consider a basis $\cB$ of $N$ and write it in the form
$$
\cB=\{x_k : k \in \N, k \ge 1\} \cup \cT
$$
where $\cT$ is infinite.

Take a finite ordered tuple $\av y$ of pairwise
distinct elements of $\cT$ of length $c-2.$
We claim that the following elements $\gamma_1,\gamma_2$ of $\iat_{c-1}(N)$ where
$$
\begin{array}{llllll}
\gamma_1: \quad & x_1      &\mapsto x_1 [x_1,\av y, x_1]             & \gamma_2: \quad &  x_1        &\mapsto x_1, \\
                & x_2      &\mapsto x_2 [x_2,x_1, \av y]\inv\quad    &                 &  x_2        &\mapsto  x_2[x_2,x_1,\av y], \\
                & x_3      &\mapsto x_3                              &                 &  x_3        &\mapsto x_3[x_3,\av y,x_1]\inv, \\
                &          &\vdots                                   &                 &             & \vdots \\
                & x_{3k+1} &\mapsto x_{3k+1} [x_{3k+1},\av y, x_1]   &                 &  x_{3k+1}   &\mapsto x_{3k+1}, \\
                & x_{3k+2} &\mapsto x_{3k+2} [x_{3k+2},x_1, \av y]\inv\quad &          &  x_{3k+2}   &\mapsto x_{3k+2}[x_{3k+2},x_1,\av y], \\
                & x_{3k+3} &\mapsto x_{3k+3}                         &                 &  x_{3k+3}   &\mapsto x_{3k+3}[x_{3k+3},\av y,x_1]\inv, \\
                &          &\vdots                                   &                 &             & \vdots \\
                & t        &\mapsto t                                &                 &  t          &\mapsto t \quad (t \in \cT).
\end{array}
$$
are commutators in the group $\ia N.$ We base our
argument on the proof of Lemma 6 of \cite{BrGu}. Let $\alpha,\beta \in \ia N$
be defined as follows:
$$
\begin{array}{llllll}
\a:\quad & x_{3k+1} &\mapsto x_{3k+1} [x_{3k+2},x_1],\quad   & \beta: \quad & x_{3k+1} &\mapsto x_{3k+1} \\
         & x_{3k+2} &\mapsto x_{3k+2},                       &              & x_{3k+2} &\mapsto  x_{3k+2} [x_{3k+1},\av y], \\
         & x_{3k+3} &\mapsto x_{3k+3},                       &              & x_{3k+3} &\mapsto x_{3k+3}, \\
         & t        &\mapsto t,                              &              & t        &\mapsto  t, \qquad (t \in \cT)
\end{array}
$$
where $k$ is an arbitrary natural number. Using the elegant idea of \cite{BrGu},
we first calculate the images of the basis elements under $\a \b$
and then under $\beta \alpha.$ Then we introduce elements $\xi_1,\xi_2$
of $\iat_{c-1}(N)$ such that $\xi_1$ takes $x_{3k+1}$ to $\gamma_1(x_{3k+1})$
$(k \in \N)$ and fixes all other elements of $\cB$ and $\xi_2$ which takes
$x_{3k+2}$ to $\gamma_1\inv(x_{3k+2})$ ($k \in \N)$ and again fixes
all other elements of $\cB.$ Next, one sees that the products
$\alpha\beta$ and $\beta\alpha$ can be ``balanced''
by placing $\xi_1,\xi_2$ into corresponding ``pans'':
$$
\xi_1 \alpha \beta = \xi_2 \beta \alpha.
$$
As $\gamma_1=\xi_1 \xi_2\inv,$ we have the desired.
A similar argument proves that $\gamma_2$ is also
a commutator in the group $\ia N.$

The product $\varepsilon=\gamma_1 \gamma_2$ is
$$
\begin{array}{lll}
\varepsilon: \quad      & x_1      &\mapsto x_1 [x_1,\av y,x_1],\\
			& x_2      &\mapsto x_2,\\
			& x_3      &\mapsto x_3[x_3,\av y,x_1]\inv, \\
			&          &\vdots \\
			& x_{3k+1} &\mapsto x_{3k+1} [x_{3k+1},\av y, x_1],\\
			& x_{3k+2} &\mapsto x_{3k+2}, \\
			& x_{3k+3} &\mapsto x_{3k+3}[x_{3k+3},\av y,x_1]\inv \\
			&          &\vdots \\
			& t        &\mapsto t \quad (t \in\cT).
\end{array}
$$
Now the set consisting of $x_i$ with $i \ne 1$ going to $x_i [x_i,\av y,x_1]$ under
$\varepsilon$ (say, the set of ``pluses''), the set consisting
of $x_i$ going to $x_i [x_i,\av y,x_1]\inv$
under $\varepsilon$ (the set of ``minuses'') and the set consisting
of $x_i$ going themselves are all moieties of $\{x_k : k \in \N, k \ge 1\}.$
Then there is a permutational automorphism $\pi$ fixing $x_1$ and
all elements of $\cT$ such that
$$
\begin{array}{lll}
\varepsilon \varepsilon^\pi: \quad & x_1 &\mapsto x_1 [x_1,\av y,x_1]^2,\\
				   & b   &\mapsto b \qquad (b \in \cB \setminus \{x_1\})
\end{array}
$$
($\pi$ must just interchange the set of ``pluses'' and the set of ``minuses''). Let $\delta_c$
denote the automorphism
$$
\begin{array}{lll}
\delta_c : \quad                   & x_1 &\mapsto x_1 [x_1,\av y,x_1],\\
				   & b   &\mapsto b \quad (b \in \cB \setminus \{x_1\})
\end{array}
$$

Fix till the end of the proof an automorphism $\Delta$ of the group $\aut N$ stabilizing
all conjugations. Due to the remark we have made in the beginning
of the proof, $\Delta$ stabilizes the automorphism $\varepsilon\varepsilon^\pi=\delta^2_c,$
an element of the commutator subgroup of the group $\iat(N).$
On the other hand, $\Delta(\delta_c)=\delta_c \eta$ for some $\eta \in \iat_{c-1}(N).$
Hence
$$
\delta_c^2=\Delta(\delta_c^2)=\delta_c^2 \eta^2,
$$
whence $\eta=\id_N.$ Thus $\delta_c$ is fixed by $\Delta.$

Recall that given a relatively
free group $F$ of infinite rank and a basis $\cX$
of $F,$ an automorphism $\s \in \aut F$ is called
{\it finitary} relative to $\cX$ if $\s$ fixes
pointwise a cofinite subset of $\cX.$

It then follows from Theorem \ref{BrGu_T,deltas} that all finitary automorphisms
in the group $K=\iat_{c-1}(N)$ (with respect to the basis $\cB$) are fixed
by $\Delta.$

Take an $\eta \in \iat_{c-1}(N).$ We apply Proposition \ref{IA_k-stabs_Are_Def}.
Observe that $\eta$ is fully determined
by the family of finitary automorphisms $\{\eta_b : b \in \cB\}$
such that $\eta_b d =d$ for all $d \in \cB \setminus \{b\}$
satisfying the conditions
\begin{equation}
\eta \in \eta_b K_{(b)}\qquad (b \in \cB).
\end{equation}
Now by Proposition \ref{IA_k-stabs_Are_Def}, $\Delta(K_{(b)})=K_{(b)},$
since $\Delta$ fixes conjugation $\tau_b,$ and
$\Delta(\eta_b)=\eta_b$  ($b \in \cB).$ Then
$\eta$ is stabilized by $\Delta.$

Write $I$ for $\ia N$ and $I_k$ for $\iat_k(N).$
We use downward induction on $k.$
Suppose that all elements of the subgroup $I_{k}$
are fixed by $\Delta$ where $k \ge 3.$ We consider an element
$$
\begin{array}{lll}
\delta_{k} : \quad             & x_1 &\mapsto x_1 [x_1,\av y,x_1],\\
				   & b   &\mapsto b \quad (b \in \cB \setminus \{x_1\})
\end{array}
$$
of $I_{k-1}$ where $\av y$ is a tuple of pairwise distinct
elements of $\cT$ of length $k-2.$ Consider the free
nilpotent group $\widetilde N=N/N_{k+1}$; let denote $\widetilde{\phantom m}$
the homomorphism $\aut N \to \aut{\widetilde N},$ determined
by the natural homomorphism $N \to \widetilde N.$ As the
image of $\delta^2_{k}$ under $\widetilde{\phantom m}$ is in
the commutator subgroup of the group $\ia{\widetilde N},$
$$
\delta^2_{k} =\varepsilon \beta
$$
where $\varepsilon$ is in $I'$ and $\beta \in I_{k}.$
By the induction hypothesis, the automorphism
in the right-hand side is fixed by $\Delta,$
and then $\delta_{k}$ is also fixed by
$\Delta.$

By applying Theorem \ref{BrGu_T,deltas} once again,
we see that all finitary automorphisms of $N$ lying
in the group $I_{k-1}$ are contained in the group
$$
\nc( I' \cap I_{k-1}, \delta_{k}) \cdot I_{k},
$$
and hence all are stabilized by $\Delta.$ We then complete
the proof as in the case of the subgroup $I_{c-1}.$
\end{proof}

\section{Automorphisms of $N$ that fix a very few primitive elements}

In this section we shall obtain
a group-theoretic characterization of the conjugacy
class of the permutational
automorphism $\pi^*$ we have used extensively in the previous
section. Recall that $\pi^* \in \aut N$ was defined as
a permutational automorphism with respect to
some basis $\cX$ such that there was exactly one element $x \in \cX$
that was fixed by $\pi^*,$ whereas there were infinitely
many infinite orbits (cycles) of $\pi^*$
on $\cY =\cX \setminus \{x\}.$ It is easy to see
that $x,x\inv$ are the only primitive elements
fixed by $\pi^*$ and conjugations $\tau_x,\tau_{x\inv}$
are the only conjugations by primitive elements
fixed by $\pi^*$ under the conjugation action (Claim \ref{Pi_star!!!}).
Hence having conjugation by a primitive element
fixed by $\pi^*$ we can choose in a unique way one of the primitive
elements it is determined by. In the next section we shall use this fact to ``build''
over $\pi^*$ the ``edifice'' of the stabilizer $\Gamma_{(x)}$
of $x$ in the group $\Gamma=\aut N.$

\begin{Lem} \label{Def_o_Pi_star}
An automorphism $\pi \in \aut N$ is a conjugate
of $\pi^*$ if and only if

{\em (a)} there is a basis set $X$ of conjugations
such that some $\tau' \in X$ is fixed by $\pi$
and $\pi$ acts on $X \setminus \{\tau'\}$ as
a permutation with empty fixed-point
set and with infinitely many infinite orbits;

{\em (b)} there is a linear order $<$ on $X$
invariant under $\pi$ such that
for every $n \ge 2$ there is a $\rho_n \in \sym X$
which also preserves $<$ and $\pi=\rho^n_n.$
\end{Lem}

\begin{proof} $(\Leftarrow).$ We prove
the sufficiency part first, partly in order
to clarify the role of the condition (b)---for there is no difficulty
in the constructing automorphisms $\rho_n$
for any conjugate of $\pi^*.$

\def\trle{\triangleleft}

Let us fix some basis $\cB$ of $N.$ We start with the following observation: suppose
there is a linear order $\triangleleft$ on $\cB$ which
is invariant under some $\s \in \sym \cB.$ We then
claim that there is a basis $\cU_\s$ of the group
$N_c,$ also invariant under $\s.$ Such a basis,
as we shall see right away, can be constructed with the use
of the well-known collecting process.

Let $F(A)$ be a free group with a basis $A$ which
is assumed to be linearly ordered.
Recall that the collecting process produces the sets $\cK_n(A)$ of basic commutators
of weight $n$ where $n \ge 1.$ The sets $\cK_n(A)$ are constructed by induction as follows:
$\cK_1(A)=A,$ and to construct the set $\cK_n(A)$ of basic commutators
of weight $n,$ we assume that the set $\bigcup_{m < n} K_m(A)$
is linearly ordered in a way that the restriction of the order on $A$
is the given one, $u < v$ for every $u \in \cK_m(A)$
and every $v \in \cK_l(A)$ with $m < l,$ and the restriction
of the order on the set $\cK_m(A)$ where $m < n$ could be any.
The set $\cK_{n}(A)$ is then formed from all commutators $[u,v]$ of weight $n$
such that
\begin{equation} \label{CollProc}
\begin{array}{l}
\text{$u \in \cK_m(A),$ $v \in \cK_l(A)$ with $m+l=n$ and $u > v,$ and, moreover,} \\
\text{if $u=[u_1,u_2],$ then $u_2 \le v.$}
\end{array}
\end{equation}
 It can be proven
that the image of the set $\cK_m(A)$ under the natural
homomorphism $F(A) \to L=F(A)/\gamma_{m+1}(F(A)),$
from $F(A)$ onto the free nilpotent group $L$ of class
$m,$ is a basis of the free abelian group $\gamma_m(L)$
(see \cite[Section 4.3]{Bakht}, \cite[Section
3.1]{HN}).  It is then rather clear that the collecting
process could be applied to any basis of a free nilpotent group
$L$ of class $c,$ and the $c$-th step will produce a basis of
the $c$-th term $L_c$ of the lower central series of $L.$

Let us apply the collecting process to $\cB$
to construct a required basis of $N_c$ invariant under $\s.$
The set $\cK_2(\cB)$ consists of the commutators $[b_i,b_j]$
where $b_i,b_j \in \cB$ such that $b_i\, \triangleright\, b_j.$ As
the order $\triangleleft$ is invariant under $\s,$
$\s \cK_2(\cB)=\cK_2(\cB),$
and, moreover, the lexicographic-like order on $\cK_2(\cB)$ (after
identifying commutators $[a,b]$ with 2-letter words $ab$) is also
invariant under $\s.$

We then proceed by induction on $m,$ constructing
the sets $\cK_m(\cB),$ where $m \le c,$ all equipped
with linear orders, namely, with the lexicographical-like orders
determined by the orders on the sets $\cK_l(\cB)$
with $l < m.$ As above, we introduce the lexicographical-like order
on $\cK_m(\cB)$ where $m \ge 3$ by identifying the basic commutators
$[u,v]=[[u_1,u_2],v]$ with 3-letter words $u_1u_2v.$
The invariance of these orders under $\s$ will
guarantee the invariance of the sets $\cK_m(\cB)$ under $\s$ (for every
commutator $[u,v] \in \cK_m(\cB),$ the elements
$\s u, \s v$ satisfy the conditions (\ref{CollProc}),
and hence $[\s u,\s v]$ is in $\cK_m(\cB)).$

Next, we deal with a relevant problem on automorphisms
of free abelian groups.

Consider a free abelian group $G$ with a basis
$$
\{f\} \cup \{e_k : k \in \Z\}.
$$
Let $H$ denote the subgroup generated by the elements $e_k.$
Consider an automorphism $\rho$ of $G$ with
$\rho(e_k)=e_{k+1}$ for all $k \in \Z$ such that
$\rho$ fixes $f$ modulo $H.$

Let $H_0$ denote the subgroup
of $H$ with the basis
\begin{equation}
\{e_k-\rho(e_k) : k \in \Z\}=\{e_k-e_{k+1} : k \in \Z\}.
\end{equation}
Clearly, $H=\str{e_0} \oplus H_0.$ Suppose
that $\rho(f) \not\equiv f \Mod{H_0}.$ As, however,
$\rho(f) \equiv f \Mod H,$ then
$$
\rho(f)-f = \a e_0 + h_0
$$
where $\a$ is a nonzero integer number and $h_0 \in H_0.$
We may assume, without loss of generality, that the index of the first
nonzero coefficient in the representation of $h_0$
as a linear combination of vectors in (\theequation) is
$\ge 0.$ Indeed, let
$$
\rho(f)-f= \a e_0 + \gamma_m( e_m-e_{m+1}) +\ldots
$$
where $\gamma_m$ is the first nonzero coefficient. We have
the desired if $m \ge 0.$ If, for instance, $m=-1,$ then
$$
\a e_0+\gamma_{-1} (e_{-1}-e_0)+\ldots =
\a (e_{-1}+e_0-e_{-1}) + \gamma_{-1} (e_{-1}-e_0)+\ldots =
\a e_{-1} +(\gamma_{-1}-\alpha) (e_{-1}-e_0)+\ldots.
$$
If $m < -1,$ a trick like that has to be
repeated several times.

So let
\begin{equation} \label{actn_in_a_free_ab}
\rho(f)=f+\a e_0+\beta_0(e_0-e_1)+\ldots+\beta_r(e_r-e_{r+1})
\end{equation}
where $\beta_k$ are integers and, as we assumed before,
$\a \ne 0.$

\begin{Claim} \label{ATerrifThing:)}
Let $n \ge 1$ and
$$
\rho^n(f)=f+g_n
$$
where $g_n \in H.$ Then $g_n$
has at least $n$ nonzero coordinates with regard
to the basis $\{e_k : k \in \Z\}$ of $H.$
\end{Claim}

\begin{proof}
Rewriting (\ref{actn_in_a_free_ab}) as
$$
\rho(f+\beta_0 e_0+\ldots+\beta_r e_r)=f+\beta_0e_0+\ldots+\beta_re_r +\a e_0
$$
one easily deduces that
$$
\rho^n(f)=f+\a e_0+\ldots+\a e_{n-1}+\beta_0(e_0-e_n)+\ldots+\beta_r(e_r-e_{r+n}),
$$
whence
$$
g_n=\a e_0+\ldots+\a e_{n-1}+\beta_0(e_0-e_n)+\ldots+\beta_r(e_r-e_{r+n}).
$$
Let first $r \le n-1.$ Then
$$
g_n=\underbrace{(\a+\beta_0)e_0+\ldots+(\a+\beta_r)e_r}_{a}+
\a e_{r+1}+\ldots+\a e_{n-1}-\underbrace{(\beta_0 e_n+\ldots+\beta_r e_{r+n})}_{b}.
$$
Suppose that exactly $k$ (out of $r+1$) coordinates of the vector $a$
are zero. This means that the vector $b$ has at least $k$ nonzero
coordinates. Thus the number of nonzero coordinates of $g_n$
is at least
$$
(r+1-k)+(n-1-r)+k=n.
$$

Now let $r \ge n.$ We have that
$$
g_n = (\a+\beta_0)e_0+\ldots+(\a+\beta_{n-1})e_{n-1}+\sum_{k=n}^r (\beta_k-\beta_{k-n}) e_k
-\beta_{r-n+1} e_{r+1}-\ldots-\beta_r e_{r+n}.
$$
Observe that the numbers
$$
r-n+1,\ldots,r
$$
form a complete set of representatives modulo $n.$ Let $m$ be an element
of this set. Suppose that $m=qn+s$ where $0 \le s < n$ and that
\begin{equation}
\beta_m=0, \beta_m-\beta_{m-n}=0, \beta_{m-n}-\beta_{m-2n}=0, \ldots, \beta_{m-(q-1)n}-\beta_{s}=0.
\end{equation}
It follows that $\beta_s=0,$ and hence the $s$-th coordinate,
that is, $(\a+\beta_s)$ is nonzero, or, to put it more generally,
there is an $l$ with $l \equiv m \Mod n$ such that
the $l$-th coordinate is nonzero. Clearly, the latter also holds,
if the condition (\theequation) is not true.
Thus $g_n$ has at least $n$ nonzero coordinates, as claimed.
\end{proof}

Let $\pi \in \aut N$ satisfy the conditions (a) and (b).
Let $\cX$ be a basis of $N$ determining $X;$ we assume
that $x$ is an element of $\cX$ with $\tau'=\tau_x.$

The order $<$ mentioned in (b)
determines the linear order, say $\trle$ of
$\cX$:
$$
y \trle z \iff \tau_y < \tau_z \quad (y,z \in \cX).
$$
Now the order $\trle$ on $\cX$ is invariant
under $\pi$ modulo $N_c.$ According to our considerations
above, there is a basis $\cU_\pi$ of the group $N_c,$
consisting of basic commutators over $\cX,$ which is also invariant
under $\pi$ (strictly speaking, the arguments above
guarantee a basis in question for some $\widetilde \pi$
congruent to $\pi$ modulo the subgroup $\iat_{c-1}(N);$
but then the actions of $\pi$ and $\widetilde \pi$ on $N_c$ are identical).
Clearly, $\pi$ acts on $\cU_\pi$ as a permutation
with empty fixed-point set and infinitely many
infinite orbits.

As $\pi \tau_x \pi\inv=\tau_x,$ we have that $\pi x=x t$ where $t \in N_c.$
We are going to show that there is a $v \in N_c$ with $t=v \pi(v\inv).$
It will follow that $\pi(xv)=xv,$ that is, that $\pi$ fixes
a primitive element as $\pi^*$ does, completing the proof.

Let $n \ge 2$ be a natural number. As $\rho_n,$
a root of $\pi$ of degree $n,$ preserves
the order on $X$ that $\pi$ does, the basis
$\cU_\pi$ must also be invariant under $\rho_n,$
and the action of $\rho_n$ on this basis must also be similar
to that one of $\pi$: with infinitely many
infinite orbits and with empty fixed-point set.

Let
$$
\cU_\pi = \bigcup_{i\in I} O_i
$$
be the partition of $\cU_\pi$ into $\rho_n$-orbits.
We are in a position to apply Claim \ref{ATerrifThing:)}---to
the free abelian group $\str{x,N_c}$ with a basis $\{x\} \cup \cU_\pi$
and to the restrictions
on this group of automorphisms under consideration.
Assume that
\begin{equation}
\rho_n x = x \prod_{i \in I} u_{i}^{\a_i} v_{i} \rho_n(v_{i}^{-1})
\end{equation}
where, of course, only finitely many terms in
the latter product are $\ne 1$ and for all $i,$ $u_{i} \in O_{i},$ $\a_i \in \Z,$
and $v_{i}$ is in the
subgroup generated by $O_{i}.$

By Claim \ref{ATerrifThing:)}, if at least one of the exponents
$\a_i$ is nonzero, then the length of $t$ with regard
to the basis $\cU_\pi$ is at least $n.$ It then follows
that there is an $n$ for which $\a_i$ in (\theequation)
are all zero, and hence $\rho_n(x)=x v \rho_n(v\inv)$
for a suitable $v \in N_c.$ But this means that
$\pi(xv)=\rho_n^n(xv) =xv.$

$(\Rightarrow).$  Consider the automorphism
$$
f(q)=q+1 \quad (q \in \Q)
$$
of $(\Q;<),$ the set $\Q$ of
rational numbers equipped with the natural order $<$.
Clearly, $f$ admits extraction of roots of any natural degree $n \ge 1$
in the group $\aut{\Q;<}$: an automorphism $g_n \in \aut{\Q;<}$ where
$$
g_n(q)=q+\frac1n \quad (q \in \Q)
$$
is such that $g^n_n=f.$

Note that $f$ acts on $\Q$ without fixed points
and has infinitely many infinite orbits.

Take a basis $\cX =\{x\} \cup \bigcup_{i \in I} \cY_i$
of $N$ where $I$ is a linearly ordered set of power $|N|$
and the elements of any set $\cY_i$ having countably infinite cardinality are indexed by rational
numbers:
$$
\cY_i = \{y_{i,q} : q \in \Q\}.
$$
We then introduce a linear order $\trle$ on $\cX$ by setting
that
$$
y_{i,q} \trle x \quad (i \in I, q \in \Q)
$$
and that
$$
y_{i,q} \trle y_{i',q'} \iff (i,q) <_{\text{\scriptsize lex}} (i',q')
$$
for all $i,i' \in I$ and $q,q' \in \Q,$
where $<_{\text{\scriptsize lex}}$ is the lexicographic order on $I \times \Q.$

We define $\pi \in \sym\cX$ via
$$
\begin{array}{l}
\pi x =x,\\
\pi y_{i,q}=y_{i,f(q)}, \quad (i \in I, q \in \Q).
\end{array}
$$
As $f \in \aut{\Q;<},$ $\pi$ the preserves $\trle.$ Clearly,
$\pi$ is a conjugate of $\pi^*.$

Let $n \ge 1$ be a natural number.
Now the automorphism $\rho_n$ of $N$ such that
$$
\begin{array}{l}
\rho_n x =x,\\
\rho_n y_{i,q}=y_{i,g_n(q)}, \quad (i \in I, q \in \Q),
\end{array}
$$
is a root of $\pi$ of degree $n$ preserving $\trle.$
Finally, both $\pi$ and $\rho_n$ preserve the order $<$ on
$X=\{\tau_z : z \in \cX\}$ where
$$
\tau_{z} < \tau_{z'} \iff z \trle z'
$$
for all $z,z' \in \cX.$
\end{proof}

We record the following simple observation about $\pi^*$
for future reference.

\begin{Claim} \label{Pi_star!!!}
Let $\cX$ be a basis of $N$ such that $\pi^* \in \sym\cX$ and
there is exactly
one element $x \in \cX$ fixed by $\pi^*,$
while the action of $\pi^*$ on $\cX \setminus \{x\}$
is one with infinitely many infinite orbits.
Then $x$ and $x\inv$ are the only primitive
elements fixed by $\pi^*$ {\em(}fixed by
$\pi^*$ up to congruence modulo any subgroup $N_k$
where $2 \le k \le c${\em).}
\end{Claim}

\begin{proof}
It is convenient to prove the statement by induction
on $c.$ When $c=1,$ $N$ is a free abelian group
and the result is obvious. Suppose that
$\pi^*(a)=a$ where $a$ is a primitive element
of $N.$ Consider the natural homomorphism $\widehat{\phantom a} : N \to N/N_c.$
This homomorphism determines the homomorphism
$\aut N \to \aut{N/N_c}$ which we denote by the
same symbol $\widehat{\phantom a}.$ Then
$\widehat{\mathstrut \pi^*}(\widehat{\mathstrut a}) =\widehat{\mathstrut a},$
whence $\widehat a =\widehat x^{\pm 1}.$
Therefore $a =x^{\pm 1}t,$ where $t=t(x,\av y) \in N_c$ and
$\av y$ is a tuple of element of $\cX \setminus \{x\}.$
It follows that
$$
t(x,(\pi^*)^k \av y)=t(x,\av y),
$$
for every $k \in \Z$ and then $t(x,\av y)=1.$
\end{proof}

\section{Reconstruction of primitive elements} \label{Rec_o_PrimEls}

The aim of this section is to show definability
of stabilizers $\Gamma_{(x)}$ of primitive elements of $N$
where $\Gamma=\aut N.$

\begin{Th} \label{Def_o_Stabs}
{\em (i)} Let $z$ be a primitive element of $N$ which
is fixed by $\pi^*.$ Then the group $\Gamma_{(z)}$ of all automorphisms of $N$
stabilizing $z$ is definable with the
parameter $\pi^*;$

{\em (ii)} the family of stabilizers of
primitive elements of $N$ is invariant
under the action of the group $\aut{\aut N}.$
\end{Th}

\begin{proof}
(i). We first note that by Claim \ref{Pi_star!!!} $\pi^*$ fixes exactly
two primitive elements that are inverses of one
another. So the element $z$ in the condition
of (i) can equivalently be replaced with $z\inv,$ because
$\Gamma_{(z)}=\Gamma_{(z\inv)}.$

\def\Zpps{Z^+(\pi^*)}

Let $\tau^*$ be conjugation by a primitive
elements stabilized by $\pi^*$ under the
conjugation action. Again, it follows from Claim \ref{Pi_star!!!}
(for the case when $k=c$) that there are exactly
two such conjugations, inverses of one another. Let $x$ be the (only) primitive
element of $N$ such that $\pi^*(x)=x$ and $\tau^*=\tau_x.$
Recall that we write $\Zpps$ for the subgroup
of the centralizer $Z(\pi^*)$ of $\pi^*$ consisting
of automorphisms of $N$ that fix $x.$ Then the subgroup
$\Zpps$ is definable with the parameter
$\pi^*,$ for $\s \in \Zpps$ if and only
if $\s$ preserves $\tau^*$ under the conjugation
action and commutes with $\pi^*.$

Let us fix a basis $X$ of $\inn N$ containing
conjugation $\tau^*$ and a basis $\cX$ of $N,$
containing $x,$ which determines $X.$
Write $\cY$ for the set $\cX \setminus \{x\}.$

We are going to prove that the stabilizer $\Pi_{(x)}$
of $x$ in the group $\Pi=\sym{\cX}$ is contained
in a definable with the parameter $\pi^*$ subgroup
of $\aut N.$ We shall then demonstrate that the subgroup
$\Pi_{(x)}$ and the subgroup $Z^+(\pi^*)$ generate
the required subgroup $\Gamma_{(x)}$:
$$
\Gamma_{(x)}=\str{\Pi_{(x)},Z^+(\pi^*)}.
$$

All actions mentioned in the following
lemma are actions on $X$ by conjugation.

\begin{Lem}
{\em (a)} Let $\psi \in \sym X$ be an involution
which fix $\tau^*=\tau_x$ and which acts on $X$ as a transposition
interchanging distinct conjugations $\tau_y=\tau_y(\psi)$ and
$\tau_z=\tau_z(\psi)$ of $X$ where $y,z \in \cY.$
Suppose that $\s$ is a power of $\pi^*$
such that $\tau_y$ and $\tau_z$ do not lie
in the same orbit of $\s$ and for every natural number $n \ge 1$
there is a $\rho_n \in \sym X$ having
the following properties:
\begin{itemize}
\item $\rho_n$ commutes with the $n$-th
power $\s^n$ of $\s$;
\item the orbits $O_n,O_n' \sle X$ of $\s^n$
containing $\tau_y$ and $\tau_z,$ respectively,
are both fixed pointwise by $\rho_n$;
\item the actions of $\s^n$ and $\rho_n$
on $X \setminus (O_n \cup O_n')$ are the same;
\item $\psi$ commutes with $\rho_n.$
\end{itemize}
Then $\psi x=x t$ where $t$ is an element
of $N_c \cap \str{x,y,z}.$

{\em (b)} Suppose that for an involution
$\psi$ satisfying the conditions of {\em (a)}
there is an involution $\rho \in \sym X,$ commuting
with a power of $\pi^*,$ which takes $\tau_y(\psi)$
to $\tau_z(\psi)$ and commutes with $\psi.$ Then $\psi x=x.$

{\em (c)} Let $\Psi \sle \sym X$ be a family of involutions
satisfying the conditions of {\em (b)} such that
for every distinct $y,z \in \cY$ there is exactly
one element of $\Psi$ taking $\tau_y$ to $\tau_z.$
Suppose that an involution $\f \in \sym X$
preserves $\tau^*$ and normalizes
$\Psi$: $\f \Psi\f =\Psi.$ Then $\f x=x.$
\end{Lem}

\begin{proof}
(a) Assume that $\psi x =x t$ where $t \in N_c$ and that
$t$ is written as a reduced word in letters of $\cX.$

We work with the case when $n=1.$ Note that
being an element of $\sym X,$ $\rho_1$ cannot
invert $\tau^*.$ Then as $\rho_1$ commutes with $\s,$
we have that $\rho_1 x=x.$ Further, $\rho_1 \psi \rho_1\inv=\psi$
implies that $x\rho_1(t)=xt,$ or $\rho_1(t)=t.$
Now $\rho_1$ acts a permutation having no fixed
points and having infinitely many infinite orbits
on $X \setminus (\{\tau_x\} \cup O_1 \cup O_1').$
Let subsets $\cO_1$ and $\cO_1'$ of $\cX$ induce
the sets $O_1$ and $O_1'.$ It then follows that the reduced word $t$
has no occurrences of letters in $\cX \setminus (\{x\} \cup\cO_1 \cup \cO_1'),$
or $t \in N_c \cap \str{\{x\} \cup \cO_1 \cup \cO_1'}.$

Let $\cO_n,\cO_n'$ be the subsets of $\cX$ determining the sets $O_n,O_n',$ respectively.
Arguing as above, we see that $t \in N_c \cap \str{\{x\}\cup \cO_n \cup \cO_n'}$
for all $n \ge 1.$
Observe that
$$
\bigcap_{n \ge 1} (\cO_n \cup \cO_n') = \{y,z\}.
$$
Hence $t=t(x,y,z).$

(b). As $\psi$ is an involution, then $\psi^2 x=x$
implies that $x t(x,y,z) t(x,z,y)=x,$
whence $t(x,z,y)=t(x,y,z)\inv.$
 Let $\rho$ satisfy the conditions of (b). Then
$$
\rho \psi \rho(x) =\psi x \iff \rho[x t(x,y,z)] = x t(x,y,z) \iff
x t(x,z,y)=xt(x,y,z).
$$
Then $t(x,y,z)^2=1,$ and $t=1,$ as claimed.

(c) Suppose that $\f x =xs$ where $s \in N_c.$ Clearly,
$\f$ sends $s$ to the inverse. Choose an arbitrary
$\psi \in \Psi.$ We then have that $\f \psi \f =\psi'$
for a suitable $\psi' \in \Psi.$ Then $\f \psi \f x=\psi' x =x.$
Therefore
$$
x =\f \psi \f (x) =\f \psi (xs)=\f(x \psi s) =x \cdot s \cdot \f \psi(s),
$$
or
$$
\f \psi(s) =s\inv \iff \psi(s)=\f(s\inv) \iff \psi(s)=s.
$$
Thus $\psi(s)=s$ for all $\psi \in \Psi.$
This easily implies that $s=1,$ for the subgroup
of $\sym X$ generated by $\Psi$ acts on $\cX \setminus \{x\}$ modulo $N_c$
$k$-transitively for every natural number $k \ge 1.$
\end{proof}

Now we are constructing a promised definable with the parameter
$\pi^*$ subgroup $S$ containing the group $\Pi_{(x)}.$
Statements of Lemma \theLem\ suggest the construction
of $S$ as the subgroup of $\aut N$ generated by all involutions
in $\sym X$ that are in the normalizers
of families involutions satisfying the conditions of (c):
\begin{equation} \label{def_o_S}
S =\str{\f \in \sym X : \f^2 =\id, \f \tau^* \f =\tau^*, \exists \Psi \text{ with (c) s.t. } \f \Psi \f=\Psi}.
\end{equation}
By Lemma \theClaim, $\s x=x$ for all $\s \in S.$ Let us
show that
$$
S \sge \Pi_{(x)} =\{\s \in \sym\cX : \s x =x\}.
$$

Indeed, consider the family $\Psi^*$ of transpositions $\psi_{y,z} \in \sym\cX$ where distinct
$y,z$ run over $\cY$ such that $\psi_{y,z}$ interchanges
$y$ with $z$ and fixes all other elements of $\cX.$

Let us check that $\Psi^*$ satisfies (c). Take any $\psi_{y,z} \in \Psi^*.$
The check of (a) for $\psi_{y,z}$ is simple: one easily
finds a power $\s$ of $\pi^*$ such $y,z$ lie in the different
orbits of $\s.$ Then the construction of a needed
$\rho_n$ for every $n \ge 1$ is trivial. Let us check (b).
Consider as a $\rho$ an element of $\sym\cX$
which fixes $x,$ and sends $y$ to $z$ while
interchanging the orbits of $\s$ on $\cX$
containing $y$ and $z$; all other orbits may, say stay
fixed pointwise under the action of $\rho.$  Hence the conditions of (b) are also
true for $\psi_{y,z}.$

Recall that any infinite symmetric group is generated by
involutions (as a bit too powerful tool one can apply,
for instance, Baer's theorem on the structure
of normal subgroups of the infinite symmetric
group $\sym\Omega$ \cite[p. 256]{DiMo}; it is an immediate corollary
of the said theorem that the normal
closure of any element of $\sym\Omega$---in particular,
any involution---whose support
is of power $|\Omega|$ is the whole of
$\sym\Omega).$

Now as any involution $\f$ of the group $\Pi_{(x)}$ normalizes $\Psi^*,$
the group $\Pi_{(x)}$ is contained in $S.$

Let us again consider a relevant problem
on free abelian groups. Let $G$
be a free abelian group with a basis
$\{a\} \cup \{b_i : i \in I\}.$ Let $H$ denote
the subgroup generated by $\{b_i : i \in I\}.$
Then an arbitrary
automorphism $\beta$ of $G$ fixing $a$ looks like
$$
\begin{array}{l}
\beta a = a\\
\beta b_i = c_i +k_i a, \quad (i \in I)
\end{array}
$$
where $c_i$ is in $H$ and $k_i \in \Z$ ($i \in I$).
Clearly, if $\alpha \in \aut G$ is defined via
$$
\begin{array}{l}
\a a = a\\
\a b_i = b_i -k_i a, \quad (i \in I),
\end{array}
$$
then $\beta \alpha \in \aut G{}_{(a),\{H\}}.$

This observation implies that
$$
\Gamma_{(x)} = \Gamma_{(x),\{\str\cY\}} \cdot A_x \cdot \iat(N)_{(x)}.
$$
where $A_x$ is the group of all automorphisms of $N$ of the
form
\begin{equation}
\begin{array}{l}
\a x =x,\\
\a y = y x^{k_y}, \quad (y \in \cY)
\end{array}
\end{equation}
where $k_y \in \Z$ for all $y \in \cY.$

We then claim that
$$
\Gamma_{(x)} =\str{\Pi_{(x)},Z^+(\pi^*)}.
$$
This will complete the proof of (i).

Indeed, the group $\ia N{}_{(x)}=\iat_1(N)_{(x)}$ is contained
in the group $\str{\Pi_{(x)},Z^+(\pi^*)}$
by Lemma \ref{I_k_in_the_LHS}. Let $\cZ$ be any moiety of $\cY.$
It is easy to see that $A_x$ equals
the product the subgroup $(A_x)_{(\cX \setminus \cZ)}$
and the conjugate subgroup of $(A_x)_{(\cX \setminus \cZ)}$
by a suitable element of $\Pi_{(x)}.$ But then,
by Lemma \ref{Populating_cZ_0},
the group $\str{\Pi_{(x)},Z^+(\pi^*)}$
contains the group $(A_x)_{(\cX \setminus \cZ_0)}$
for a suitable moiety $\cZ_0$ of $\cY.$

The same lemma implies that the group
$\str{\Pi_{(x)},Z^+(\pi^*)}$ contains
the subgroup $\Gamma_{\str{\cZ_0}}$ where
$\cZ_0$ is a moiety of $\cY.$
As the author showed in \cite{To:JLM_Bergman}, given
any partition $\cX=\cX_1 \cup \cX_2 \cup \cT$ of $\cX$ into
pairwise disjoint moieties, we have that the group $\Gamma=\aut N$
is generated by the stabilizers
\begin{equation}
\Gamma_{(\cX_1),\{\str{\cX_2,\cT}\}} \text{ and }
\Gamma_{(\cX_2),\{\str{\cX_1,\cT}\}}.
\end{equation}
(Theorem 2.5 in \cite{To:JLM_Bergman}).
Clearly, the stabilizers in (\theequation)
are conjugate by a suitable permutational
automorphism of $N$ with respect to the basis $\cX.$
Thus $\aut N$ is generated by some $\rho \in \sym\cX$
and the stabilizer $\Gamma_{\str{\cX_2,\cT}}=\Gamma_{(\cX_1),\{\str{\cX_2,T}\}}.$

By applying these considerations to the situation
at hands, we see that the subgroup
$\Gamma_{(x),\{\str{\cY}\}} \cong \aut{\str \cY}$ is generated
by an element of $\Pi_{(x)}$
and the group $\Gamma_{\str{\cZ_0}} \sle \str{\Pi_{(x)},Z^+(\pi^*)}.$

(ii) Let $\Delta \in \aut{\aut N}$ and
$\Sigma \sle \aut N$ be the stabilizer
of a primitive element of $N.$ Then
$\Delta(\Sigma)$ is also a stabilizer
of a primitive element of $N$ by
Lemma \ref{Def_o_Pi_star} and (i).
\end{proof}

 \section{Stabilizing everything}

The family of all pairs $(\tau,\Sigma)$ where $\tau$
and $\Sigma$ satisfy the condition
\begin{equation}
\begin{array}{l}
\text{``$\tau$ is conjugation determined by a primitive element of $N$ and} \\
\text{$\Sigma \sle \aut N$ is the stabilizer of a primitive element such that} \\
\text{every element of $\Sigma$ commutes with $\tau$''}
\end{array}
\end{equation}
\noindent is a definable object over $\aut N$
by Theorem \ref{Def_o_conjs}, Lemma \ref{Def_o_Pi_star} and Theorem \ref{Def_o_Stabs}.
We shall call pairs with (\theequation) {\it primitive} pairs.

Clearly, every primitive pair $(\tau,\Sigma)$ uniquely determines some
primitive element $p(\tau,\Sigma)$ of $N.$ Indeed, a conjugation $\tau$
determines the coset $x N_c$ where $\tau=\tau_x,$
consisting of primitive elements and
we then choose the only element $xs$ where $s \in N_c$
which is stabilized by every member of $\Sigma.$

Now the conjugation action
$$
\s * (\tau,\Sigma)=(\s \tau \s\inv, \s \Sigma \s\inv)
$$
of the group $\aut N$ on the set of all primitive pairs
is equivalent to the action of the group $\aut N$
on the set of all primitive elements of $N.$

We are going to show
that the multiplication of ``independent'' primitive
elements is definable.
So let us take two pairs $(\tau_1,\Sigma_1)$ and $(\tau_2,\Sigma_2)$
with (\theequation) which determine primitive
elements $x$ and $y$ of $N,$ respectively. To explain
that $x$ and $y$ are independent, we require that
there be a basis set of conjugations into
which both conjugations $\tau_1,\tau_2$ can be included.

\begin{Lem} \label{Def_o_Mult_o_Ind_Prims}
The primitive pair $(\tau_1\tau_2,\Sigma_{1,2})$ associated
with the product $xy$ of the elements $x=p(\tau_1,\Sigma_1)$ and $y=p(\tau_2,\Sigma_2)$ is
definable with the parameters $(\tau_1,\Sigma_1), (\tau_2,\Sigma_2).$
\end{Lem}

\begin{proof}
Let $X$ be any basis set of conjugations that contains
both $\tau_1$ and $\tau_2.$ For every $\tau \in X \setminus \{\tau_1,\tau_2\}$
we choose a pair $(\tau,\Sigma_\tau)$ with (\theequation) determining
an element $z_\tau$, thereby defining a basis
$$
\cX=\{x,y\} \cup \{z_\tau : \tau \in X \setminus \{\tau_1,\tau_2\} \}
$$
of $N.$

Consider a pair $(\tau_1\tau_2,\Sigma)$ satisfying (\theequation).
It is clear that this pair determines a primitive element $p$ of the
form $p=xyt$ where $t \in N_c.$ Let $w(*_1,*_2)$ be a group word in
two letters over an alphabet having no common elements with $N$
(or a term of the language of group theory having
two variables) such that $t=w(x,y)$ where $w(x,y)$ is the result of replacing
letters/variables of $w$ to $x$ and $y$: $*_1 \to x,$ $*_2 \to y.$

Observe that given two distinct elements $z_1,z_2$
of $\cX,$ the element $r=z_1 z_2 w(z_1,z_2)$ is definable
over $\cX$ and $p,$ for $r$ is the result of action
on $p=xyw(x,y)$ by a permutational automorphism
with respect to $\cX$ which takes $x$ to $z_1$
and $y$ to $z_2.$ For instance, the element $r=yxw(y,x)$
must be described as follows: first, we are describing the
uniquely determined $\pi \in \aut N$ such that
$$
\begin{array}{l}
\pi * (\tau_1,\Sigma_1)=(\tau_2,\Sigma_2), \\
\pi * (\tau_2,\Sigma_2)=(\tau_1,\Sigma_1), \\
\pi * (\tau,\Sigma_\tau)=(\tau,\Sigma_\tau),\quad (\tau \in X \setminus \{\tau_1,\tau_2\}),\\
\end{array}
$$
and then set the pair $R$ corresponding to $r$
as $R=\pi * (\tau_1\tau_2,\Sigma).$ On similar occasions below, we will just describe actions
on primitive elements themselves, skipping translations into
the language of action on primitive pairs.

Our first condition, we are going to impose on the pair $(\tau_1\tau_2,\Sigma),$
will imply that $w(x,y)=w(y,x).$
To achieve that we require the following automorphism $\alpha$ defined via
$$
\begin{array}{ll}
\alpha: \quad & x \mapsto x \\
	      & yx w(y,x) \mapsto  xy w(x,y), \\
	      & zx w(z,x) \mapsto xz w(x,z), \qquad z \in \cX \setminus \{x,y\} \\
\end{array}
$$
be equal to $\tau_1=\tau_x.$ It then follows that
$$
\begin{array}{l}
\alpha(x)=x \\
\alpha(y)=xyx\inv w(x,y)w(y,x)\inv, \\
\alpha(z)=xzx\inv w(x,z)w(z,x)\inv, \qquad z \in \cX \setminus \{x,y\}.
\end{array}
$$
whence
\begin{equation}
w(x,y)=w(y,x),
\end{equation}
as claimed. Clearly, as $x,y$ are members of some
basis of $N,$ then $w(a,b)=w(b,a)$ for every
$a,b \in N$ (apply to the both parts of (\theequation)
an endomorphism of $N$ taking $x$ to $a$ and $y$ to $b.$)

Fix a $z$ from $\cX \setminus \{x,y\}$ and
write $q$ for $xzw(x,z).$ Consider the automorphisms
$$
\begin{array}{ll}
U_1: & y \mapsto yzw(y,z), \\
     & a \mapsto a, \qquad (a \in \cX \setminus \{y\})
\end{array}
$$
and
$$
\begin{array}{ll}
U_2: & x \mapsto xyw(x,y), \\
     & a \mapsto a, \qquad (a \in \cX \setminus \{x\});
\end{array}
$$
We then require $U_1(p)$ and $U_2(q)$ be equal:
$$
U_1(p) = xyz w(y,z) w(x,yz) = U_2(q)= xyz w(x,y) w(xy,z).
$$
Therefore
$$
w(x,y) w(xy,z)=w(y,z) w(x,yz).
$$
Fix a natural number $k$ and apply the endomorphism of $N$
taking $y$ to $x^k$ and fixing all other elements
of $\cX$ to the both parts of the last equation. It then
follows that
\begin{equation}
w(x^{k+1},z)=w(x^k,z) w(x,x^k z).
\end{equation}
We apply a permutational automorphism with respect
to $\cX$ interchanging $x$ and $z$ to the both parts of
(\theequation):
$$
w(z^{k+1},x)=w(z^k,x) w(z,z^k x).
$$
Due to ``symmetricity'' of $w$ we then have that
$$
w(x,z^{k+1})=w(x,z^k) w(z^k x,z).
$$
Multiplying the latter equation with (\theequation) part by
part, we arrive at
$$
w(x^{k+1},z)w(x,z^{k+1})=w(x^k,z) w(x,z^k) \cdot w(x^k z,x) w(z^k x,z).
$$
Then by (\ref{glue_em_to_kill_em}),
$$
w(x,z)^{(k+1)^c} = w(x,z)^{k^c} w(xz^k,z) w(zx^k,x),
$$
or
$$
w(x,z)^{(k+1)^c-k^c} =w(xz^k,z) w(zx^k,x).
$$
Let $L_k=\{s^k : s \in N_c\}$ be the
subgroup of all $k$-th powers of elements of the (free abelian) group
$N_c.$ As $(k+1)^c -k^c \equiv 1 \Mod k$ and as, for instance,
$w(xz^k,z) \equiv w(x,z) \Mod{L_k},$ we get that
$$
w(x,z) \equiv 1\Mod{L_k}
$$
for every natural number $k.$ Therefore $w(x,z)=1,$
which means that the element $xyt$ determined by
the pair $(\tau_1 \tau_2, \Sigma)$ is $xy,$ as desired.

{\sc Remark.} The definability of $xy$ over $\{x,y\}$ is much
more simple when the nilpotency class $c$ of $N$
is {\it even.} Indeed, extend $\{x,y\}$
to a basis $\cX$ of $N$ and consider
the automorphism $\theta \in \aut N$ that inverts all elements of $\cX.$
Then the element $xyt$ where $t \in N_c$ is inverted
by the automorphism $\tau_x \theta$ if and only
if $t=1.$ This condition can be easily translated
into the language of the action of $\aut N$ on the
primitive pairs.

\end{proof}

\begin{Th}
The group $\aut N$ is complete.
\end{Th}

\begin{proof}
According to Proposition \ref{Centreless}, the
group $\aut N$ is centreless.

Let $\Delta \in \aut{\aut N}.$ By Proposition
\ref{StabilizingConjugations} there is a $\s_0 \in \aut N$
such that $\Delta_1=T_{\s_0} \circ \Delta$ stabilizes
all elements of the group $\inn N$
and, according to the same Proposition,
\begin{equation}
\Delta_1(\s) \equiv \s \Mod{\iat_{c-1}(N)} \quad (\s \in \aut N).
\end{equation}

Take a basis set $X$ of conjugations, and
then for each $\tau \in X$ choose the stabilizer
$\Sigma_\tau$ of a primitive element,
forming a primitive pair $(\tau,\Sigma_\tau).$ Let $P(X)$ denote the set
$$
\{(\tau,\Sigma_\tau) : \tau \in X\}
$$
and let $\cX$ be a basis of $N$ determined
by the pairs in $P(X).$ Suppose
that a pair $(\tau,\Sigma_\tau)$ determines
the primitive element $x_\tau.$

Take a $\tau \in X.$
The image of the pair $(\tau,\Sigma_\tau)$ under $\Delta_1$
is a pair $(\tau,\Sigma_\tau^*)$ which determines
a primitive element $x_\tau t_\tau$ where $t_\tau \in N_c.$
Indeed, let $\pi,$ a conjugate of $\pi^*,$
belong to $\Sigma_\tau;$ then $\pi \tau \pi\inv=\tau$
by the definition of a primitive pair. By (\theequation),
$$
\Delta_1(\pi) \tau \Delta_1(\pi)\inv =\tau.
$$
On the other hand, $\Delta_1(\pi)$ must be
a conjugate of $\pi^*$ (Lemma \ref{Def_o_Pi_star}) and
hence $\Delta(\pi)$ fixes a uniquely determined element in the coset $x_\tau N_c,$
the above mentioned element $x_\tau t_\tau.$

Consider the automorphism $\a \in \iat_{c-1}(N)$ which
takes $x_\tau$ to $x_\tau t_\tau$ ($\tau \in X).$
Then $\Delta_2 = T_{\a\inv} \circ \Delta_1$ stabilizes
all pairs $(\tau,\Sigma_\tau)$ where $\tau$ runs over
$X.$

\begin{Lem} \label{Inn_then_IA}
$\Delta_2$ {\em (}and hence already $\Delta_1${\em)} stabilizes all elements of the
group $\ia N.$
\end{Lem}

\begin{proof}
Let $x,y$ be two distinct elements of $\cX.$ Consider
the IA-automorphism $K_{12}$ (the standard notation) defined via
$$
\begin{array}{l}
K_{12} (xy)=yx \iff K_{12}(x)=yxy\inv,\\
K_{12} z =z, \quad (z \in \cX \setminus \{x\}).
\end{array}
$$
Then as the pairs associated with the elements
$xy$ and $yx$ are definable over $P(X)$
and $P(X)$ is stabilized by $\Delta_2,$
$K_{12}$ is stabilized by $\Delta_2.$

According to the classical results of Nielsen and Magnus,
the normal closure of an analogue of $K_{12}$ in
the automorphism group $\aut{F_n},$ where $F_n$ ($n \ge 2$)
is a finitely generated free group, is the
group $\ia{F_n}$ \cite[Section 3.5]{MKS}. Also, the homomorphism $\aut{F_n} \to
\aut{F_n/\gamma_3 (F_n)},$ determined
by the natural homomorphism $F_n \to F_n/\gamma_3(F_n),$
is surjective \cite{Andrea} (equivalently, one says
that the automorphism groups of finitely
generated free nilpotent groups of class
two are tame). It follows that in the
automorphism group of any infinitely
generated free nilpotent group of class two all finitary
IA-automorphisms are contained in the
normal closure of any analogue of $K_{12}.$
Therefore finitary IA-automorphisms of $N$
are all contained in the group $\nc(K_{12}) \cdot \iat_2(N).$

By Proposition \ref{StabilizingConjugations}, all elements of $\nc(K_{12})$ are stabilized
by $\Delta_2$ and the group $\iat_2(N)$ is stabilized
by $\Delta_2$ by Proposition \ref{StabilizingIA3}.
Thus all finitary IA-automorphisms of $N$
are stabilized by $\Delta_2,$ and by using
the fact that subgroups $\ia N_{(p)}$
where $p$ is any primitive element of $N$
are invariant under $\Delta_2$ (Proposition \ref{IA_k-stabs_Are_Def}), we can
complete the proof as in the proof
of Proposition \ref{StabilizingIA3}.
\end{proof}

Swan (see \cite[Section 2]{BurnsPi}) found rather simply-constructed
generating sets of the automorphism groups
of infinitely generated free abelian groups.
Let $A$ be an infinitely generated free abelian group
and let $\cB$ be a basis of $A.$ An
automorphism $\theta$ of $A$ is called
{\it $\cB$-block-unitriangular} if there is a moiety
$\cC$ of $\cB$ such that $\theta$ fixes $\cC$ pointwise and
$$
\theta d \equiv d~(\mathop{\rm mod} \str{\cC}).
$$
for all $d \in \cB \setminus \cC.$ Then the group $\aut A$ is generated by
all $\cB$-block-unitriangular automorphisms.

Let $\cY$ be a moiety of $\cX.$
Consider then an automorphism $\beta$ of $N$ which induces
a block-unitriangular automorphism in the abelianization
of $N$:
\begin{equation}
\begin{array}{ll}
\beta z = z \prod_{y \in \cY} y^{m(y,z)}\quad &(z \in \cX \setminus \cY),\\
\beta y = y \quad                             &(y \in \cY)
\end{array}
\end{equation}
where in any (formal) product $\prod_{y \in \cY} y^{m(y,z)}$ only
finitely many integer exponents $m(y,z)$ are non-zero.

By Swan's theorem, the group $\aut N$ is generated
by the automorphisms of the form (\theequation) and
by the elements of the group $\ia N.$ It therefore
remains to prove that $\Delta_2$ stabilizes
all automorphisms of the form (\theequation).

We use Lemma \ref{Def_o_Mult_o_Ind_Prims}. Any element
of the form $z \prod_{y \in \cY} y^{m(y,z)}$ can be
obtained from elements of $\cX$ as a result of a number of successive applications of the operation
``multiplication of independent primitives''
like in the following example:
$$
z \to z \cdot y_1^{-1} \to zy_1^{-1} \cdot y_1^{-1} \to
zy_1^{-2} \cdot y_2 \to \ldots \to zy_1^{-2} y_2^2 \cdot y_2 \to \ldots
$$
Therefore the image under $\beta \in \aut N$ with (\theequation)
of any basis pair $(\tau,\Sigma_\tau)$
that determines an element of $\cX \setminus \cY$
is a primitive pair definable with the parameters
from the set of basis pairs $P(X).$ As $\Delta_2$ fixes
all members of $P(X),$ it fixes $\beta$ and the result follows:
$\Delta_2$ is then the trivial, and $\Delta$ is an inner
automorphism of the group $\aut N.$
\end{proof}

{\sc Acknowledgements.}
The author would like to thank Oleg
Belegradek for helpful discussions and to
thank the referee for the constructive
comments on the manuscript.

\end{document}